\documentclass[leqno,a4paper,12pt]{article}
\usepackage{amsmath,amsthm}

\newtheorem{Thm}{\indent Theorem}[section]
\newtheorem{Prop}[Thm]{\indent Proposition}
\newtheorem{Lem}[Thm]{\indent Lemma}
\newtheorem{Cor}[Thm]{\indent Corollary}

\theoremstyle{definition}
\newtheorem{Def}[Thm]{\indent Definition}
\newtheorem{Rem}[Thm]{\indent Remark}

\newtheorem{Ass}[Thm]{\indent Assumption}

\def\qed{{\hskip0pt\unskip\unskip\nobreak\hfil\penalty50
          \hskip1em\hbox{}\nobreak\hfil
          {\bf q.e.d.}%
          \parfillskip=0pt\finalhyphendemerits=0
          \par}\medskip}

\newenvironment{Proof}
               {{\it Proof.}\quad}
               {\qed}

\newenvironment{Proofof}[1]
               {{\it Proof of #1.}\quad}
               {\qed}

\newcommand{\Prime}{\kern3\fontdimen1\font$'$\kern-7\fontdimen1\font}

\long\def\forget#1{}

\long\def\beginSIDEREMARK#1\endSIDEREMARK
    {{\par\bigskip\advance\leftskip by 2cm
                  \advance\rightskip by -2cm\noindent
      {\bf Our own side remark:} #1
      \par\bigskip\noindent}}

\long\def\beginFORGET#1\endFORGET{#1}
\long\def\beginFORGET#1\endFORGET{}

\def\?{\ ???\ \immediate\write16{}%
\immediate\write16{Warning: There was still a question mark . . . }%
\immediate\write16{}}

\usepackage{amsmath}

\usepackage{exscale}
\usepackage{amssymb}

\newcommand{\BA}{{\mathbb{A}}}

\newcommand{\BC}{{\mathbb{C}}}

\newcommand{\BF}{{\mathbb{F}}}
\newcommand{\BG}{{\mathbb{G}}}

\newcommand{\BN}{{\mathbb{N}}}

\newcommand{\BQ}{{\mathbb{Q}}}
\newcommand{\BR}{{\mathbb{R}}}
\newcommand{\BS}{{\mathbb{S}}}

\newcommand{\BZ}{{\mathbb{Z}}}

\newcommand{\Fc}{{\mathfrak{c}}}

\newcommand{\Fn}{{\mathfrak{n}}}
\newcommand{\Fo}{{\mathfrak{o}}}
\newcommand{\Fp}{{\mathfrak{p}}}

\newcommand{\FA}{{\mathfrak{A}}}
\newcommand{\FB}{{\mathfrak{B}}}

\newcommand{\FH}{{\mathfrak{H}}}

\newcommand{\FL}{{\mathfrak{L}}}

\newcommand{\FS}{{\mathfrak{S}}}

\newcommand{\FX}{{\mathfrak{X}}}
\newcommand{\FY}{{\mathfrak{Y}}}
\newcommand{\FZ}{{\mathfrak{Z}}}

\newcommand{\CA}{{\cal A}}

\newcommand{\CC}{{\cal C}}
\newcommand{\CD}{{\cal D}}

\newcommand{\CG}{{\cal G}}
\newcommand{\CH}{{\cal H}}

\newcommand{\CL}{{\cal L}}

\newcommand{\CN}{{\cal N}}

\newcommand{\CS}{{\cal S}}

\newcommand{\CV}{{\cal V}}
\newcommand{\CW}{{\cal W}}

\newcommand{\ch}{\mathop{\rm CH}\nolimits}
\newcommand{\diag}{\mathop{\rm diag}\nolimits}
\newcommand{\Spec}{\mathop{{\bf Spec}}\nolimits}

\newcommand{\rad}{\mathop{{\rm rad}}\nolimits}

\newcommand{\End}{\mathop{\rm End}\nolimits}

\newcommand{\GL}{\mathop{\rm GL}\nolimits}
\newcommand{\Gm}{\mathop{\BG_m}\nolimits}
\newcommand{\GFm}{\mathop{\BG_{m,F}}\nolimits}
\newcommand{\Gr}{\mathop{\rm Gr}\nolimits}

\newcommand{\Hom}{\mathop{\rm Hom}\nolimits}

\newcommand{\Rep}{\mathop{\rm Rep}\nolimits}

\newcommand{\SL}{\mathop{\rm SL}\nolimits}

\newcommand{\loccit}{[loc.$\;$cit.]}

\def\tei{\, | \,}
\def\halb{\frac{1}{2}}

\def\id{{\rm id}}

\newbox\mybox
\def\arrover#1{\mathrel{
       \setbox\mybox=\hbox spread 1.4em{\hfil$\scriptstyle#1$\hfil}
       \vbox{\offinterlineskip\copy\mybox
             \hbox to\wd\mybox{\rightarrowfill}}}}
\def\larrover#1{\mathrel{
       \setbox\mybox=\hbox spread 1.4em{\hfil$\scriptstyle#1$\hfil}
       \vbox{\offinterlineskip\copy\mybox
             \hbox to\wd\mybox{\leftarrowfill}}}}

\def\ontoover#1{\mathrel{
       \setbox\mybox=\hbox spread 1.4em{\hfil$\scriptstyle#1$\hfil}
       \vbox{\offinterlineskip\copy\mybox
             \hbox to\wd\mybox{\rightarrowfill\hskip-2.8mm
                               $\rightarrow$}}}}
\def\leftontoover#1{\mathrel{
       \setbox\mybox=\hbox spread 1.4em{\hfil$\scriptstyle#1$\hfil}
       \vbox{\offinterlineskip\copy\mybox
             \hbox to\wd\mybox{$\leftarrow$\hskip-2.8mm
                               \leftarrowfill}}}}
\def\longto{\longrightarrow}
\def\into{\hookrightarrow}

\def\longonto{\ontoover{\ }}
\def\isoto{\arrover{\sim}}

\def\longinto{\lhook\joinrel\longrightarrow}

\usepackage[curve,matrix,arrow,cmtip]{xy}
\NoComputerModernTips

\def\myxymessage{\def\messagetext
   {Here an xy-pic diagram was omitted to speed up compilation . . . }
   \immediate\write16{\messagetext}
   \hbox{\bf \messagetext}}
\def\filxymatrix#1{\myxymessage}
\def\filxyarray#1{\myxymessage}

\newdir^{ (}{{}*!/-3pt/\dir^{(}}
\newdir_{ (}{{}*!/-3pt/\dir_{(}}
\newdir^{ )}{{}*!/+3pt/\dir^{)}}
\newdir_{ )}{{}*!/+3pt/\dir_{)}}

\def\rscript#1{\hbox to 0pt{$\scriptstyle#1$\hss}}

\let\oldbullet\bullet
\def\bullet{{\mathchoice{\oldbullet}%
                        {\oldbullet}%
                        {\scriptscriptstyle\oldbullet}%
                        {\oldbullet}}}

\newcommand{\argdot}{{\;\bullet\;}}

\newcommand{\bA}{\mathop{\overline{A}}\nolimits}

\newcommand{\bX}{\mathop{\overline{X}}\nolimits}

\newcommand{\ua}{\mathop{\underline{\alpha}}\nolimits}

\newcommand{\CHeffM}{\mathop{CHM^{eff}(k)}\nolimits}
\newcommand{\CHM}{\mathop{CHM(k)}\nolimits}
\newcommand{\CHSM}{\mathop{CHM^s(S)}\nolimits}

\newcommand{\CHeffFM}{\mathop{\CHeffM_F}\nolimits}
\newcommand{\CHFM}{\mathop{\CHM_F}\nolimits}
\newcommand{\CHFbM}{\mathop{CHM(\bar{k})_F}\nolimits}
\newcommand{\CHFQAbM}{\mathop{CHM^{Ab}(F)_\BQ}\nolimits}
\newcommand{\CHFAbM}{\mathop{CHM^{Ab}(k)_F}\nolimits}
\newcommand{\CHFbAbM}{\mathop{CHM^{Ab}(\bar{k})_F}\nolimits}
\newcommand{\CHFSM}{\mathop{CHM^s(S)_F}\nolimits}
\newcommand{\DeffgM}{\mathop{DM^{eff}_{gm}(k)}\nolimits}
\newcommand{\DeffFgM}{\mathop{\DeffgM_F}\nolimits}

\newcommand{\DgM}{\mathop{DM_{gm}(k)}\nolimits}
\newcommand{\DFgM}{\mathop{\DgM_F}\nolimits}
\newcommand{\DFFgM}{\mathop{DM_{gm}(F)_F}\nolimits}
\newcommand{\DFQgM}{\mathop{DM_{gm}(F)_\BQ}\nolimits}
\newcommand{\DAbgM}{\mathop{DM_{gm}^{Ab}(k)}\nolimits}
\newcommand{\DFAbgM}{\mathop{DM_{gm}^{Ab}(k)_F}\nolimits}
\newcommand{\DFFAbgM}{\mathop{DM_{gm}^{Ab}(F)_F}\nolimits}
\newcommand{\DFQAbgM}{\mathop{DM_{gm}^{Ab}(F)_\BQ}\nolimits}

\newcommand{\ReF}{\mathop{{\rm Res}_{F/\BQ}}\nolimits}
\newcommand{\Mgm}{\mathop{M_{gm}}\nolimits}

\newcommand{\Mcgm}{\mathop{M_{gm}^c}\nolimits}

\newcommand{\dMgm}{\mathop{\partial M_{gm}}\nolimits}

\begin{document}

%

\hfuzz=3pt
\overfullrule=10pt                   


\setlength{\abovedisplayskip}{6.0pt plus 3.0pt}
\setlength{\belowdisplayskip}{6.0pt plus 3.0pt}
\setlength{\abovedisplayshortskip}{6.0pt plus 3.0pt}
\setlength{\belowdisplayshortskip}{6.0pt plus 3.0pt}

\setlength{\baselineskip}{13.0pt}
\setlength{\lineskip}{0.0pt}
\setlength{\lineskiplimit}{0.0pt}

%
%

\title{On the interior motive of certain Shimura varieties:
the case of Picard surfaces
\forget{
\footnotemark
\footnotetext{To appear in ....}
}
}
\author{\footnotesize by\\ \\
\mbox{\hskip-2cm
\begin{minipage}{6cm} \begin{center} \begin{tabular}{c}
J\"org Wildeshaus \footnote{
Partially supported by the \emph{Agence Nationale de la
Recherche}, project ``R\'egulateurs
et formules explicites''. }\\[0.2cm]
\footnotesize Universit\'e Paris 13\\[-3pt]
\footnotesize Sorbonne Paris Cit\'e \\[-3pt]
\footnotesize LAGA, CNRS (UMR~7539)\\[-3pt]
\footnotesize F-93430 Villetaneuse\\[-3pt]
\footnotesize France\\
{\footnotesize \tt wildesh@math.univ-paris13.fr}
\end{tabular} \end{center} \end{minipage}
\hskip-2cm}
\\[2.5cm]
}
\date{February 9, 2015}
\maketitle
\begin{abstract}
\noindent
The purpose of this article is
to construct a Hecke-equivariant Chow motive whose realizations equal
interior (or intersection) cohomology of Picard 
surfaces with regular algebraic coefficients. 
As a consequence, we are able to define Grothendieck motives
for Picard modular forms.  \\

\noindent Keywords: Picard surfaces, weight structures, 
boundary motive, interior motive.

\end{abstract}


\bigskip
\bigskip
\bigskip

\noindent {\footnotesize Math.\ Subj.\ Class.\ (2010) numbers: 
14G35
(11F32, 11F55, 14C25, 14F25, 19E15, 19F27).
}

\eject

\tableofcontents

\bigskip


%
%

\setcounter{section}{-1}
\section{Introduction}
\label{Intro}



The purpose of this paper is the construction and analysis of the 
\emph{interior motive} of Kuga--Sato families over Picard surfaces. 
It gives a second example of Shimura varieties where the use 
of the formalism of \emph{weight structures} \cite{Bo} proves
to be successful for dealing with a problem, for which explicit
geometrical methods seem inefficient, the first example being that
of Hilbert--Blumenthal varieties of arbitrary dimension \cite{W5}. \\

Let $k$ be a field of characteristic zero, 
and denote by $\DgM$ the triangulated
category of \emph{geometrical motives} \cite{V}. 
(We shall actually consider the $F$-linear variant $\DFgM$ of $\DgM$, 
for suitable coefficients $F$ of characteristic zero.
For the purpose of this introduction, let us agree not to distinguish 
notation between $\DgM$ and $\DFgM$.) 
For a smooth scheme $X$ over $k$,
denote by $\Mgm (X)$ and $\Mcgm (X)$
the \emph{motive} of $X$ and its \emph{motive with compact support}, 
respectively \cite{V}. \emph{Chow motives} form a full sub-category
of $\DgM$; indeed they are identified with the category
of objects which are \emph{pure of weight zero} with respect to
the motivic weight structure \cite{Bo}.
Both functors $\Mgm$ and $\Mcgm$ agree on smooth and proper 
$k$-schemes $X$, and yield the Chow motive of $X$. In general,
that is, for schemes $X$ which are only supposed smooth,
they are related by a natural transformation $\iota: \Mgm \to \Mcgm$. 
The morphism $\iota_X: \Mgm (X) \to \Mcgm (X)$ clearly factors 
through a Chow motive; indeed,
the motive of \emph{any} smooth compactification of $X$ will provide such
a factorization. \\

The r\^ole of the interior motive \cite{W4} is to give a canonical, and in fact
minimal such factorization. \\
 
The \emph{boundary motive} $\dMgm(X)$ of $X$ \cite{W2}
is a canonical choice of (the shift by [-1] of) a cone of $\iota_X$; 
indeed, it fits into a canonical exact triangle
\[
(\ast) \quad\quad
\dMgm(X) \longto \Mgm(X) \stackrel{\iota_X}{\longto}
 \Mcgm(X) \longto \dMgm(X)[1] \; .
\] 
In order to formulate a reasonable condition sufficient to guarantee
the existence of the interior motive, we need to
assume in addition that an idempotent endomorphism
$e$ of the exact triangle $(\ast)$ is given. We thus get a direct factor
\[
\dMgm(X)^e \longto \Mgm(X)^e \stackrel{\iota_X^e}{\longto} 
\Mcgm(X)^e \longto \dMgm(X)^e[1] \; .
\] 
Now the abstract yoga of triangulated categories shows that isomorphism
classes of factorizations $\Mgm(X)^e \to M_0 \to \Mcgm(X)^e$
of $\iota_X^e$ correspond bijectively to isomorphism classes 
of morphisms $\dMgm(X)^e \to N$ (define $N$ as the shift by $[-1]$ of a
cone of $M_0 \to \Mcgm(X)^e$). \emph{A priori},
it thus appears plausible that the
existence of a minimal factorization $M_0$, which is pure of weight zero,
should be related to a condition on $\dMgm(X)^e$. \\

Indeed, a sufficient such condition was identified in \cite{W4}. It says
that the object $\dMgm(X)^e$ \emph{avoids weights $-1$ and $0$}
(in practice, this condition is never satisfied for the whole of 
$\dMgm(X)$, whence the need to consider direct factors).  \\

Technically speaking, the present paper deals with the verification
of that hypothesis in the context of Picard surfaces. 
Our task is considerably simplified by a recent result of Ancona's \cite{Anc},
which implies that the analogue of the hypothesis is satisfied
for the Hodge structure on \emph{boundary cohomology}. \\

Our main task thus consists in proving what could be named
the principle of \emph{weight conservativity}: the absence of certain weights
in the \emph{realization} $R(M)$ of a motive $M$ implies the absence of these
weights in $M$ itself. For \emph{Artin--Tate motives} over
a number field, weight conservativity is known \cite{W8}, and relatively straightforward 
to establish thanks to the presence
of a $t$-structure. However, the boundary motive of a Kuga--Sato
family over a Picard surface is not of Artin--Tate type; in fact,
motives of certain elliptic curves (with complex multiplication) do occur.
While in the general situation, a $t$-structure
is expected, its existence has not yet been established ---
not even on the sub-category of $\DgM$ generated by 
the motive of a single elliptic curve over a number field! \\

Our first main result, Theorem~\ref{1C} establishes weight conservativity
on the triangulated sub-category $\DAbgM$ of $\DgM$ generated
by \emph{Chow motives of Abelian type} over a field $k$ which can
be embedded into $\BC$. Its proof consists of several intermediate steps,
each of which may be considered to be of independent interest: (1)~first,
we show that the realization functor $R$ is conservative on the category
of Chow motives of Abelian type (Theorem~\ref{1I}), (2)~then,  
we extend conservativity to the whole of $\DAbgM$ (Theorem~\ref{1L}), 
(3)~finally, we show that conservativity actually implies Theorem~\ref{1C}. 
The proof of step~(1) uses a result from
\cite{AK}: $R$ is conservative if its source is
\emph{semi-primary} and pseudo-Abelian, if $R$ itself is
\emph{radicial}, and if the only object mapped to zero under
$R$ is the zero object. This is where our definition enters:
Chow motives of Abelian varieties are \emph{finite dimensional}
\cite{Ki}, hence according to one of the main results from
\cite{AK}, they generate a semi-primary category. The main result
from \cite{L}, proved for Abelian varieties, can be
reformulated to say that $R$ is radicial
on Chow motives of Abelian type. For steps~(2) and (3), 
we make a systematic use of the 
\emph{minimal weight filtration}, whose theory was developed in \cite{W9}. \\ 

In Section~\ref{1}, the application of Theorem~\ref{1C} to the boundary motive of 
a smooth scheme $X$ over $k$ is made explicit (Corollary~\ref{1E}). 
In order to do so, we identify
its realizations with boundary cohomology of $X$ (Proposition~\ref{1Ea}). \\   

Section~\ref{2} then contains 
our second main result, Theorem~\ref{2Main}, which can be seen
as the motivic analogue of the main result of \cite{Anc}. 
Here, the smooth scheme
$X$ is an $N$-th power of the universal Abelian threefold over a smooth
Picard surface $S$ associated to a quadratic imaginary number
field $F$. Denote by $G$ the reductive group 
underlying the Shimura variety $S$, and let $V$ be 
an irreducible algebraic representation of $G$. According to
the main result from \cite{Anc1} (valid for arbitrary Shimura varieties
of $PEL$-type), $N$ can be chosen such that the \emph{relative Chow
motive} $h(X/S)$ admits an idempotent $e$ with the following 
property: the cohomological
Hodge theoretic realization of $h(X/S)^e$ coincides with the variation
of Hodge structure on $S$ induced by $V$. 
Theorem~\ref{2Main} then implies that the boundary motive
$\dMgm(X)$ lies in $DM_{gm}^{Ab}(F)$; furthermore, its direct factor
$\dMgm(X)^e$ avoids weights $-1$ and $0$ if and only if 
the representation $V$ is of regular highest weight.
Therefore, the interior motive exists. 
We list its principal pro\-perties, using the main results
from \cite[Sect.~4]{W4}.
First (Corollary~\ref{2G}), we get precise statements on the
weights occurring in the motive $\Mgm(X)^e$ and the motive with compact
support $\Mcgm(X)^e$. Second (Corollary~\ref{2H}),
the interior motive is Hecke-equivariant. 
Corollary~\ref{2H} appears
particularly interesting, given the problem of non-existence of
equivariant smooth compactifications of $X$ (for $N \ge 1$).
Third (Corollary~\ref{2I}), the 
interior motive occurs canonically as a direct factor 
of the (Chow) motive of any smooth compactification of $X$. \\

Section~\ref{3} is devoted to the verification 
of the criterion from Theorem~\ref{1C}.
First, we need to show that in this case, the boundary motive
is indeed of Abelian type (Theorem~\ref{3A}). This is done using a 
smooth \emph{toroidal compactification}. We use \emph{co-localization}
for the boundary motive \cite{W2}, in order to reduce 
to showing the statement for the contribution of any of the strata.
The latter was identified in the general context of \emph{mixed Shimura
varieties} \cite{W3}. For Kuga--Sato families over 
Picard surfaces, \loccit \ shows 
in particular that these contributions
are indeed all of Abelian type over $F$. 
All that then remains to do is to cite the main result from \cite{Anc}. \\

In the final Section~\ref{4}, we give the necessary ingredients to
perform the construction of the Grothendieck motive associated to
a (Picard) automorphic form (Definition~\ref{4E}). 
This is the analogue for Picard surfaces
of the main result from \cite{Sc}. \\

Let us stress that 
our present approach \emph{via} weight structures
appears ne\-cessary because contrary
to the situation considered in \cite{Sc}, Hecke-equivariant
compactifications of the involved Kuga--Sato families are not known 
(and possibly not even expected) to exist. 
One may reasonably expect this approach to generalize. We refer
in particular to \cite{C} for the study of Picard varieties
of higher dimension. \\  

I wish to thank G.~Ancona, G.~Harder, A.~Mokrane,
V.~Pilloni, S.~Rozensztajn and J.~Tilouine 
for useful discussions and comments. 
Special thanks go to the referee for a detailed reading
of the first version of this article; her or his comments and suggestions
contributed considerably to improve its organization. \\

{\bf Notation and conventions}: For a perfect field $k$, 
we denote by $Sch/k$ the category of separated schemes of finite 
type over $k$, and by $Sm/k \subset Sch/k$
the full sub-category of objects which are smooth over $k$. 
As far as motives are concerned,
the notation of this paper is that of \cite{W2,W3,W4}, which in turn follows
that of \cite{V}. We refer to \cite[Sect.~1]{W2} for a concise
review of this notation, and of the de\-fi\-nition of the triangulated 
categories $\DeffgM$ and $\DgM$ of (effective) geometrical
motives over $k$. Let $F$ be a commutative $\BQ$-algebra.
The notation $\DeffFgM$ and $\DFgM$ stands 
for the $F$-linear analogues of these triangulated categories
defined in \cite[Sect.~16.2.4
and Sect.~17.1.3]{A}. 
Similarly, let us denote by $\CHeffM$ and $\CHM$ the categories 
opposite to the categories of (effective) Chow motives over $k$, and by
$\CHeffFM$ and $\CHFM$ the pseudo-Abelian
completion of the category $\CHeffM \otimes_\BZ F$ and 
$\CHM \otimes_\BZ F$, respectively. 
Using \cite[Cor.~4.2.6]{V}, we canonically identify 
$\CHeffFM$ and $\CHFM$ with
a full additive sub-category of $\DeffFgM$ and $\DFgM \, $, respectively. 
For $S \in Sm/k$, denote by
$\CHSM$ the category opposite to the category of \emph{smooth
Chow motives over $S$} \cite[Sect.~1.3., 1.6]{DM},
and by $\CHFSM$ the pseudo-Abelian
completion of the category $\CHSM \otimes_\BZ F$. 


\bigskip

%
%

\section{Motives of Abelian type}
\label{0}



Let us start by defining the categories of motives for which we are
able to establish conservativity and weight conservativity
of realizations. The notations
are those fixed at the end of the Introduction; in particular, $k$ denotes
a perfect field, and $F$ a commutative $\BQ$-algebra.

\begin{Def} \label{1B}
Denote by $\bar{k}$ the algebraic closure of $k$. \\[0.1cm] 
(a)~Define the category of \emph{Chow motives of Abelian type over $\bar{k}$}
as the strict full dense additive tensor sub-category $\CHFbAbM$ of $\CHFbM$
generated by the following:
\begin{enumerate}
\item[(1)] shifts of Tate motives $\BZ(m)[2m]$, for $m \in \BZ$, 
\item[(2)] Chow motives of Abelian varieties over $\bar{k}$. 
\end{enumerate}
(b)~Define the category of \emph{Chow motives of Abelian type over $k$}
as the (strict) full (dense additive tensor) sub-category $\CHFAbM$ of $\CHFM$ 
of Chow motives whose base change to $\bar{k}$ lies in 
$\CHFbAbM \subset \CHFbM$. \\[0.1cm] 
(c)~Define the category $\DFAbgM$ as the (strict) full triangulated 
sub-category of $\DFgM$ generated by $\CHFAbM$.
\end{Def}

Note that both $\CHFbAbM$ and $\CHFAbM$ are rigid.
The following observation will turn out to be vital.

\begin{Prop} \label{1J}
The \emph{motivic weight structure} on $\DFgM$ (see \cite[Sect.~6]{Bo}) 
induces a weight structure
on $\DFAbgM$. More precisely, we have
\[
DM_{gm}^{Ab}(k)_{F,w \le 0} = \DFAbgM \cap DM_{gm}(k)_{F,w \le 0} 
\]
and 
\[
DM_{gm}^{Ab}(k)_{F,w \ge 0} = \DFAbgM \cap DM_{gm}(k)_{F,w \ge 0} \; ,
\]
and the heart $DM_{gm}^{Ab}(k)_{F,w = 0}$ equals $\CHFAbM$.
\end{Prop}

\begin{Proof}
The category $\DFAbgM$ is generated by $\CHFAbM$ as a triangulated
category. 
The latter is contained in $\CHFM$, therefore it
is \emph{negative}, meaning that
\[
\Hom_{\DFAbgM}(M_1,M_2[i]) = \Hom_{\DFgM}(M_1,M_2[i]) = 0
\]
for any two objects $M_1, M_2$ of $\CHFAbM$, and any integer $i > 0$
\cite[Cor.~4.2.6]{V}.
Therefore,
\cite[Thm.~4.3.2~II~1]{Bo} can be applied to ensure the existence
of a weight structure $v$ on $\DFAbgM$, uniquely characterized by the
property of containing $\CHFAbM$ in its heart. Furthermore
\cite[Thm.~4.3.2~II~2]{Bo}, the heart of $v$
is equal to the category of retracts of $\CHFAbM$. But by definition,
$\CHFAbM$ is pseudo-Abelian, hence $DM_{gm}^{Ab}(k)_{F,v = 0} = \CHFAbM$.
In particular, 
\[
DM_{gm}^{Ab}(k)_{F,v = 0} \subset DM_{gm}(k)_{F,w = 0} \; .
\]
The existence of weight filtrations for the weight structure 
$v$ then formally implies that
\[
DM_{gm}^{Ab}(k)_{F,v \le 0} \subset DM_{gm}(k)_{F,w \le 0} \; ,
\]
and that
\[
DM_{gm}^{Ab}(k)_{F,v \ge 0} \subset DM_{gm}(k)_{F,w \ge 0} \; .
\]
We leave it to the reader to prove from this (cmp.~\cite[Lemma~1.3.8]{Bo})
that in fact
\[
DM_{gm}^{Ab}(k)_{F,v \le 0} = DM_{gm}^{Ab}(k) \cap DM_{gm}(k)_{F,w \le 0} 
\]
and
\[
DM_{gm}^{Ab}(k)_{F,v \ge 0} = DM_{gm}^{Ab}(k) \cap DM_{gm}(k)_{F,w \ge 0} \; . 
\]
\end{Proof}

\begin{Cor} \label{1Ja}
The category $\DFAbgM$ is pseudo-Abelian.
\end{Cor}

\begin{Proof}
The weight structure $w$ on $\DFAbgM$ is \emph{bounded} in the sense that
its heart $\CHFAbM$ generates the triangulated category $\DFAbgM$.
Furthermore, the category $\CHFAbM$ is pseudo-Abelian. 

Our claim thus follows from \cite[Lemma~5.2.1]{Bo}.
\end{Proof}

Next, let us consider realizations
(\cite[Sect.~2.3 and Corrigendum]{Hu}; see \cite[Sect.~1.5]{DG}
for a simplification of this approach). We  
shall concentrate on two
realizations (Theorem~\ref{1C} then formally generalizes to
any of the other realizations ``with weights'' considered in \cite{Hu}):
\begin{enumerate}
\item[(i)] when $k$ is embeddable into $\BC$, the Hodge theoretic realization
\[
R_\sigma : \DFgM \longto D
\]
associated to a fixed embedding $\sigma$ of $k$ into
the field $\BC$ of complex numbers. Here, $D$ is the bounded derived category 
of mixed graded-polarizable $\BQ$-Hodge structures 
\cite[Def.~3.9, Lemma~3.11]{Be}, tensored with $F$,
\item[(ii)] the $\ell$-adic realization
\[
R_{\ell} : \DFgM \longto D
\]
for a prime $\ell$ not dividing $char(k)$. 
Here, $D$ is the bounded ``derived category'' of
constructible $\BQ_\ell$-sheaves on $\Spec (k)$
\cite[Sect.~6]{E}, tensored with $F$.
\end{enumerate}

Choose and fix one of these two, denote it by $R$, and
recall that it is a contravariant tensor functor mapping
the pure Tate motive $\BZ(m)$ to the pure Hodge structure $\BQ(-m)$
(when $R=R_\sigma$) and to the pure $\BQ_\ell$-sheaf $\BQ_\ell(-m)$
(when $R=R_\ell$), respectively \cite[Thm.~2.3.3]{Hu}.
The category $D$ 
is equipped with a $t$-structure;
write $D^{t=0}$ for its heart, and 
$H^n : D \to D^{t=0}$, $n \in \BZ$, for the cohomology functors. \\

In order to analyze the functor $R$,
let us start by recalling some termino\-logy from category theory.
Following \cite[Sect.~1.1]{AK}, we refer to a category $\FA$ as 
an \emph{$F$-category} ($F =$ our fixed commutative $\BQ$-algebra), 
if the morphisms in $\FA$ are equipped with the
structure of an $F$-module, and if composition of morphisms
in $\FA$ is $F$-linear.
Similarly, a functor $T: \FA \to \FB$ is an \emph{$F$-functor} if it
is $F$-linear on morphisms. 

\begin{Def}[{\cite{Ke}}]
Let $\FA$ be an $F$-category. 
The \emph{radical} of $\FA$ is the ideal $\rad_\FA$
which associates to each pair of objects $A,B$ of $\FA$ the subset
\[
\rad_\FA (A,B) := \{ f \in \Hom_\FA (A,B) \; , \; 
  \forall \, g \in \Hom_\FA (B,A) \; , \; \id_A - gf \; \text{invertible} \}
\]
of $\Hom_\FA (A,B)$.
\end{Def}

In \cite[D\'ef.~1.4.1]{AK}, the radical is referred to as the 
\emph{Kelly radical} of $\FA$.
It can be checked that $\rad_\FA$ is indeed a two-sided ideal of $\FA$
in the sense of \cite[Sect.~1.3]{AK}, \emph{i.e.}, for each pair of objects 
$A,B$, $\rad_\FA (A,B)$ is an $F$-submodule of $\Hom_\FA (A,B)$, and for each
pair of morphisms $h: A' \to A$ and $k: B \to B'$ in $\FA$, 
\[
k \rad_\FA (A,B) h \subset \rad_\FA (A',B') \; .
\]
 
\begin{Def}[{\cite[D\'ef.~1.4.6]{AK}}] 
An $F$-functor $T: \FA \to \FB$ between $F$-categories
is said to be \emph{radicial} if $T(\rad_\FA) \subset \rad_\FB$.
\end{Def}

\begin{Def}[{\cite[D\'ef.~2.3.1]{AK}}] 
An $F$-category $\FA$ is called \emph{semi-primary} if \\[0.1cm]
(1)~for all objects $A$ of $\FA$, the radical $\rad_\FA (A,A)$ is nilpotent,
\emph{i.e.}, there exists a positive integer $N$ such that 
\[
\rad_\FA (A,A)^N = 0 \subset \End_\FA (A) \; ,
\]
(2)~the $F$-linear quotient category 
\[
\bA := \FA / \rad_\FA
\]
is semi-simple. 
\end{Def}

We then have the following result.

\begin{Prop} [{\cite[Prop.~2.3.4~f)]{AK}}] \label{1F}
Assume $F$ to be a finite direct product of fields of characteristic zero.
Let $T: \FA \to \FB$ an $F$-functor of $F$-categories. Assume the following:
\begin{enumerate}
\item[(1)] the category $\FA$ is semi-primary and pseudo-Abelian, 
\item[(2)] the functor is radicial,
\item[(3)] whenever $A$ is an object of $\FA$, and $T(A) = 0$, then $A = 0$. 
\end{enumerate}
Then $T$ is conservative.
\end{Prop}

The first step of our analysis consists in applying
Proposition~\ref{1F} to the restriction of the functor
\[
H^* \! R := (H^n \! R)_{n \in \BZ} : \DFAbgM \to \Gr_{\BZ} D^{t=0} 
\]
to $\CHFAbM \subset \DFAbgM$. Here,
we denote by $\Gr_{\BZ} D^{t=0}$ the $\BZ$-graded category associated
to the heart $D^{t=0}$ of the $t$-structure on the target $D$ of the
realization $R$. Let us check the hypotheses from Proposition~\ref{1F}.

\begin{Prop} \label{1G}
Assume $F$ to be a finite direct product of fields of characteristic zero.
All motives in $\CHFAbM$ are \emph{finite dimensional} \cite[Def.~3.7]{Ki}.
The $F$-category $\CHFAbM$ is semi-primary (and pseudo-Abel\-ian). 
\end{Prop}  

\begin{Proof}
According to \cite[pp.~54--55]{O'S2}, finite dimensionality can be checked
after base change to the algebraic closure of $k$.
Thus, the first claim follows from our definition, and from 
\cite[Thm.~(3.3.1)]{Ku}.

The second claim follows from the first, and from 
\cite[Thm.~9.2.2]{AK} (see also \cite[Lemma~4.1]{O'S1}).
\end{Proof}

This shows hypothesis~\ref{1F}~(1).
By \cite[Cor.~7.3]{Ki}, whenever
$M$ is an object of $\CHFAbM$, 
and $H^n \! R(M) = 0$ for all $n \in \BZ$, then $M = 0$,
whence hypothesis~\ref{1F}~(3).
It remains to check hypothesis~\ref{1F}~(2). \\

For the rest of the section, we assume $k$ to be a field embeddable into $\BC$,
and $F$ a finite direct product of fields of characteristic zero.
The following is a reformulation of the main result from \cite{L}.

\begin{Thm} \label{1H}
The restriction of the functor
\[
H^* \! R : \DFAbgM \to \Gr_{\BZ} D^{t=0} 
\]
to $\CHFAbM$ maps the ra\-di\-cal $\rad_{\CHFAbM}$ to zero.
In particular, it is radicial. 
\end{Thm}

\begin{Proof}
By comparison, the statement for the $\ell$-adic realization
follows from the one for the Hodge theoretic realization,
hence we may assume that $R$ equals the latter.
Base change from to $k$ to $\BC$ is a radicial
functor on finite dimensional Chow motives (\cite[Thm.~9.2.2]{AK}, \cite[p.~55]{O'S2}; see also \cite[Cor.~6.2]{W9}).
We may therefore work over $\BC$. 

Let $B$ be a proper, smooth variety over $\BC$, $r = 1,2$,
such that $\Mgm(B)$ is finite dimensional. 
The radical
\[
\rad_{CHM(\BC)_F} \bigl( \Mgm(B) , \Mgm(B) \bigr)
\]
then coincides with the maximal tensor ideal 
\[
\CN \bigl( \Mgm(B) , \Mgm(B) \bigr)
\]
\cite[Thm.~9.2.2]{AK}. It thus
consists of the classes of numerically trivial cycles on 
$B \times_{\BC} B$ \cite[Ex.~7.1.2]{AK}. 
But according to \cite[Thm.~4]{L}, numerical and
homological equivalence coincide on Abelian varieties
over $\BC$. In other words, 
\[
H^* \! R \bigl( \rad_{CHM(\BC)_F} \bigl( \Mgm(B) , \Mgm(B) \bigr) \bigr) 
= 0 
\]  
if $B$ is an Abelian variety over $\BC$.
We leave it to the reader to generalize this statement, and show that
\[
H^* \! R \bigl( \rad_{CHM(\BC)_F} \bigl( \Mgm(B)(m)[2m] , 
                                            \Mgm(B)(m)[2m] \bigr) \bigr) 
= 0 
\]   
for any Abelian variety $B$ over $\BC$ and any $m \in \BZ$. 

Our claim then follows from \cite[Lemme~1.3.2]{AK}.
\end{Proof}

\begin{Thm} \label{1I}
Assume $k$ to be a field embeddable into $\BC$,
and $F$ a finite direct product of fields of characteristic zero.
The restriction of the functor
\[
H^* \! R : \DFAbgM \to \Gr_{\BZ} D^{t=0} 
\]
to $\CHFAbM$ is conservative.
\end{Thm}

\begin{Proof}
All the hypotheses of Proposition~\ref{1F} are indeed verified.
\end{Proof}

The second step of our analysis aims at  
conservativity of $H^* \! R$ on the whole of $\DFAbgM$.  
We shall need the following technical result.

\begin{Lem} \label{1K}
Let $M$ and $N$ be objects of $\DFAbgM$, and assume that $M$ is of weights
$\ge 0$, and that $N$ is of weights $\le 0$. Then the radical
\[
\rad_{\DFAbgM} ( M , N )
\]
is contained in the kernel of $H^* \! R$.
\end{Lem}

\begin{Proof}
By comparison, we may assume that $R = R_\sigma$ equals the Hodge theoretic
realization.  
Let $f \in \rad_{\DFAbgM} ( M , N )$, 
$n \in \BZ$, and consider the morphism
\[
H^n \! R(f) : H^n \! R(N) \longto H^n \! R(M) \; .
\]
It is a morphism of graded-polarizable mixed $\BQ$-Hodge structures. 
Its source $H^n \! R(N)$
is of weights at least $n$, and its 
target  $H^n \! R(M)$, of weights at most $n$.
Furthermore, the morphism $H^n \! R(f)$ 
is strict with respect to the weight filtrations. Therefore,
it is sufficient to show that the morphism 
\[
\Gr^W_n H^n \! R(f) : \Gr^W_n H^n \! R(N) \longto \Gr^W_n H^n \! R(M)
\]
induced by $H^n \! R(f)$ is zero.  

Applying Proposition~\ref{1J}, we choose weight filtrations
\[
M_0 \stackrel{\iota_0}{\longto} M \longto M_{\ge 1} \longto M_0[1]
\]
and
\[
N_{\le -1} \longto N \stackrel{\pi_0}{\longto} N_0 \longto N_{\le -1}[1]
\]
of $M$ and $N$, with $M_0, N_0 \in \CHFAbM$, 
$M_{\ge 1} \in DM_{gm}^{Ab}(k)_{F,w \ge 1}$ and
$N_{\le -1} \in DM_{gm}^{Ab}(k)_{F,w \le -1}$.
The composition $\pi_0 \circ f \circ \iota_0$ belongs to
\[
\rad_{\DFAbgM} ( M_0 , N_0 ) = \rad_{\CHFAbM} ( M_0 , N_0 ) \; ;
\]
therefore (Theorem~\ref{1H}) it belongs to the kernel of $H^n \! R$:
\[
H^n \! R (\pi_0 \circ f \circ \iota_0) = 0 \; .
\]
But $H^n \! R(\pi_0): H^n \! R(N_0) \to H^n \! R(N)$ factors through an epimorphism 
\[
H^n \! R(N_0) \longonto \Gr^W_n H^n \! R(N) \; ,
\]
while $H^n \! R (\iota_0): H^n \! R(M) \to H^n \! R(M_0)$ factors through a monomorphism
\[
\Gr^W_n H^n \! R(M) \longinto H^n \! R(M_0) \; ,
\]   
and the composition
\[
H^n \! R(N_0) \longonto \Gr^W_n H^n \! R(N) 
\stackrel{\Gr^W_n H^n \! R(f)}{\longto} \Gr^W_n H^n \! R(M) \longinto H^n \! R(M_0)
\]
equals $H^n \! R (\pi_0 \circ f \circ \iota_0)$ \cite[Prop.~2.1.2~2]{Bo}.
This composition being zero means that $\Gr^W_n H^n \! R(f)$ is zero.
\end{Proof}

\begin{Thm} \label{1L}
Assume $k$ to be a field embeddable into $\BC$,
and $F$ a finite direct product of fields of characteristic zero.
Then the functor
\[
H^* \! R : \DFAbgM \to \Gr_{\BZ} D^{t=0} 
\]
is conservative.
\end{Thm}

\begin{Proof}
Given that $\DFAbgM$ is triangulated, and that $H^* \! R$ is cohomological,
it suffices to show that whenever $M \in \DFAbgM$, and $H^* \! R(M) = 0$,
then $M = 0$. Recall that the heart $\CHFAbM$ of the weight structure
on $\DFAbgM$ is semi-primary (Proposition~\ref{1G}). Take the \emph{minimal
weight filtration} 
\[
M_{\le 0} \longto M \longto M_{\ge 1} \stackrel{\delta}{\longto} M_{\le 0}[1]
\]
of $M$ \cite[Thm.~2.2]{W9}
($M_{\le 0} \in DM_{gm}^{Ab}(k)_{F,w \le 0}$, 
$M_{\ge 1} \in DM_{gm}^{Ab}(k)_{F,w \ge 1}$), \emph{i.e.}, 
the weight filtration of $M$ whose
isomorphism class is determined by the condition 
\[
\delta \in \rad_{\DFAbgM} (M_{\ge 1},M_{\le 0}[1]) \; .
\]
Applying Lemma~\ref{1K} to $\delta[-1]$, we see that 
$H^* \! R(\delta) = 0$. In other words, the Bockstein morphisms
in the long exact cohomology sequence for $R(M)$ are all zero. 
Now $H^* \! R(M) = 0$, therefore we see that $H^* \! R(M_{\le 0}) = 0$
and $H^* \! R(M_{\ge 1}) = 0$. Since the motivic weight structure
is bounded, we are thus reduced to the case where $M$ is pure, say of weight
$m$. But then $M[-m]$ is in the heart of $w$, \emph{i.e.} 
(Proposition~\ref{1J}), in $\CHFAbM$. Our claim thus follows from 
Proposition~\ref{1G}, and from \cite[Cor.~7.3]{Ki}.
\end{Proof}

Here is the main result of this section; it generalizes \cite[Thm.~3.11]{W8}.

\begin{Thm} \label{1C}
Let $k$ be a field embeddable into $\BC$, $F$ a finite direct product of fields 
of characteristic zero, and $R$ one of the two realizations
considered above (Hodge theoretic or $\ell$-adic). 
Then $R$ respects and
detects the weight structure $w$ on $\DFAbgM$. More precisely,
let $M \in \DFAbgM$, and $\alpha \le \beta$ two integers. \\[0.1cm]
(a)~$M$ lies in the heart $\CHFAbM$ of $w$ if and only if 
the $n$-th cohomology object $H^n \! R(M) \in D^{t=0}$ of $R(M)$ 
is pure of weight $n$, 
for all $n \in \BZ$. \\[0.1cm]
(b)~$M$ lies in $DM_{gm}^{Ab}(k)_{F,w \le \alpha}$ if and only if 
$H^n \! R(M)$ is of weights $\ge n - \alpha$, 
for all $n \in \BZ$. \\[0.1cm] 
(c)~$M$ lies in $DM_{gm}^{Ab}(k)_{F,w \ge \beta}$ if and only if 
$H^n \! R(M)$ is of weights $\le n - \beta$, 
for all $n \in \BZ$. \\[0.1cm]
(d)~$M$ is without weights $\alpha,\alpha+1,\ldots,\beta$ if and only if 
$H^n \! R(M)$ is without weights $n - \beta,\ldots,n - (\alpha+1),n - \alpha$, 
for all $n \in \BZ$.  
\end{Thm}

Here, absence of weights is  
in the sense of \cite[Def.~1.10]{W4}. In the situation of Theorem~\ref{1C}~(d),
it means that there is an exact triangle
\[ 
M_{\le \alpha - 1} \longto M \longto M_{\ge \beta + 1} 
\longto M_{\le \alpha - 1}[1]
\]
in $\DeffgM_F$, with 
$M_{\le \alpha - 1} \in DM_{gm}^{Ab}(k)_{F,w \le \alpha - 1}$,
$M_{\ge \beta + 1} \in DM_{gm}^{Ab}(k)_{F,w \ge \beta + 1}$. 

\medskip

\begin{Proofof}{Theorem~\ref{1C}}
The motivic weight structure is bounded, and the ``only if'' parts of
statements~(a)--(d) are true. It therefore suffices to prove the
``if'' part of statement~(d).

Consider the minimal weight filtrations of $M$, concentrated at weight
$\alpha$ and at weight $\beta+1$, respectively \cite[Thm.~2.2]{W9}: 
\[
M_{\le \alpha-1} \longto M \longto M_{\ge \alpha} 
\stackrel{\delta_\alpha}{\longto} M_{\le \alpha-1}[1] \; ,
\]
\[
M_{\le \beta} \longto M \longto M_{\ge \beta+1} 
\stackrel{\delta_{\beta+1}}{\longto} M_{\le \beta}[1] \; .
\]
Here, $M_{\le \alpha-1} \in DM_{gm}^{Ab}(k)_{F,w \le \alpha-1}$ \emph{etc.}, 
and both $\delta_\alpha$ and $\delta_{\beta+1}$ belong to the radical.
By orthogonality, the identity on $M$ extends to a morphism of exact
triangles
\[
\vcenter{\xymatrix@R-10pt{
        M_{\le \alpha-1} \ar[d]_{m} \ar[r] & M \ar@{=}[d] \ar[r] &
        M_{\ge \alpha} \ar[d] \ar[r]^{\delta_\alpha} & M_{\le \alpha-1}[1] \ar[d]_{m[1]}\\
        M_{\le \beta} \ar[r] & M \ar[r] &
        M_{\ge \beta+1} \ar[r]^{\delta_{\beta+1}} & M_{\le \beta}[1]
\\}}
\] 
By Lemma~\ref{1K}, both $H^* \! R(\delta_\alpha)$ and 
$H^* \! R(\delta_{\beta+1})$ are zero. Thus, the above morphism of exact
triangles induces a morphism of exact sequences
\[
\vcenter{\xymatrix@R-10pt{
        0 \ar[r] &
        H^* \! R(M_{\le \alpha-1}) \ar@{^{ (}->}[d]_{H^* \! R(m)} \ar[r] & 
        H^* \! R(M) \ar@{=}[d] \ar[r] &
        H^* \! R(M_{\ge \alpha}) \ar@{->>}[d] \ar[r] & 0 \\
        0 \ar[r] &
        H^* \! R(M_{\le \beta}) \ar[r] & 
        H^* \! R(M) \ar[r] &
        H^* \! R(M_{\ge \beta+1}) \ar[r] & 0
\\}}
\] 
Our hypothesis on weights avoided in $H^* \! R(M)$ implies that
the monomorphism $H^* \! R(m)$ is in fact an isomorphism.
But then (Theorem~\ref{1L}) so is $m$ itself. This yields a weight filtration
\[
M_{\le \alpha-1} \longto M \longto M_{\ge \beta+1} 
\longto M_{\le \alpha-1}[1] 
\]
of $M$ avoiding weights $\alpha,\alpha+1,\ldots,\beta$.
\end{Proofof}

\begin{Cor} \label{1Cb}
Let $k$ be a field embeddable into $\BC$,
and $F$ a finite direct product of fields 
of characteristic zero.
Denote by $DMAT(k)_F$ the triangulated category of \emph{Artin--Tate motives 
over $k$} \cite[Def.~1.3]{W8}. Then $R$ respects and
detects the weight structure on $DMAT(k)_F$.
\end{Cor}

Note that \cite[Thm.~3.11]{W8} only concerns the case where
$k$ is an algebraic number field.

\medskip

\begin{Proofof}{Corollary~\ref{1Cb}}
Indeed, $DMAT(k)_F$ is a full sub-category of $\DFAbgM$.
\end{Proofof}

\begin{Cor} \label{1Ca}
Let $k$ be a field embeddable into $\BC$, $F$ a finite direct product of fields 
of characteristic zero, and $\alpha \le \beta$ be two integers. 
Assume that $N \in \DeffFgM$ is a successive
extension of objects $M$ of $\DeffFgM$,
each satisfying one of the following properties.
\begin{enumerate}
\item[(i)] $M$ is without weights $\alpha,\alpha+1,\ldots,\beta$.
\item[(ii)] $M$ lies in $\DFAbgM$, and 
the cohomology object $H^n \! R(M)$ of its image $R(M)$ under $R$
is without weights $n - \beta,\ldots,n - (\alpha+1),n - \alpha$, 
for all $n \in \BZ$. 
\end{enumerate} 
Then $N$ is without weights $\alpha,\alpha+1,\ldots,\beta$.
\end{Cor}

\begin{Proof}
Apply \cite[Prop.~2.11]{W8} and Theorem~\ref{1C}~(d).
\end{Proof}


\bigskip

%
%

\section{A criterion on the existence of the interior motive}
\label{1}



Let $k$ be a field embeddable into $\BC$,
and $F$ a finite direct product of fields of characteristic zero. \\

Fix $X \in Sm/k$, and consider the exact triangle
\[
(\ast) \quad\quad
\dMgm(X) \longto \Mgm(X) \longto \Mcgm(X) \longto \dMgm(X)[1]
\]
in $\DeffgM$ \cite[Prop.~2.2]{W2}. Fix 
an idempotent endomorphism $e$ of the image of the exact triangle
in the category $\DeffgM_F$. 
Denote by $\Mgm(X)^e$, $\Mcgm(X)^e$ and $\dMgm(X)^e$ the images of $e$
on $\Mgm(X)$, $\Mcgm(X)$ and $\dMgm(X)$, respectively, considered as
objects of $\DeffgM_F$. 
Recall the following assumption.

\begin{Ass}[{\cite[Asp.~4.2]{W4}}] \label{1A}
The object $\dMgm(X)^e$ 
is without weights $-1$ and $0$.
\end{Ass}

In order to apply the results from \cite[Sect.~4]{W4},
allowing in particular to construct the \emph{interior motive},
one needs to verify Assumption~\ref{1A}. The results obtained in Section~\ref{0}
yield the following criterion.

\begin{Thm} \label{1D}
Let $k$ be a field embeddable into $\BC$,
and $F$ a finite direct product of fields of characteristic zero.
Let $X \in Sm/k$, and $\alpha \le \beta$ two integers such that $\alpha \le -1$ 
and $\beta \ge 0$.
Assume that $\dMgm(X)^e$ is a successive
extension of objects $M$ of $\DeffFgM$,
each satisfying one of the following properties.
\begin{enumerate}
\item[(i)] $M$ is without weights $\alpha,\alpha+1,\ldots,\beta$.
\item[(ii)] $M$ lies in $\DFAbgM$, and 
the cohomology object $H^n \! R(M)$ of its image $R(M)$ under $R$
is without weights $n - \beta,\ldots,n - (\alpha+1),n - \alpha$, 
for all $n \in \BZ$. 
\end{enumerate} 
Then Assumption~\ref{1A} holds. 
\end{Thm}

\begin{Proof}
This is Corollary~\ref{1Ca} for $N = \dMgm(X)^e$.
\end{Proof}

Let $X \in Sm/k$. Recall from \cite[Def.~4.1~(a)]{W4} that
$c(X,X)$ contains a canonical sub-algebra $c_{1,2}(X,X)$
(of ``bi-finite correspondences'')
acting on the exact triangle 
\[
(\ast) \quad\quad
\dMgm(X) \longto \Mgm(X) \longto \Mcgm(X) \longto \dMgm(X)[1] \; .
\]
Denote by $\bar{c}_{1,2}(X,X)$ 
the quotient of $c_{1,2}(X,X)$ by the kernel of this action. 
The algebra $\bar{c}_{1,2}(X,X) \otimes_\BZ F$ is a canonical source
of idempotent endomorphisms $e$ of $(\ast)$, and it is for such choices
that the results in \cite{W4} were formulated. However, they remain
valid in the present, more general context. \\

Assuming $e \in \bar{c}_{1,2}(X,X) \otimes_\BZ F$,
the group $H^n \! R \bigl( \dMgm(X)^e \bigr)$
equals the $e$-part of the \emph{boundary cohomology} of $X$
with respect to the \emph{natural} action (\emph{via} cycles) of the algebra
$\bar{c}_{1,2}(X,X) \otimes_\BZ F$ \cite[Prop.~2.5]{W5}.
We shall need a variant of this result. 
In order to explain it, recall that boundary cohomology is 
defined \emph{via} a compactification
$j: X \into \bX$: writing $i: \partial \! \bX \into \bX$
for the complementary immersion,
one defines $\partial H^n (\bullet)$ as cohomology of $\partial \! \bX$
with coefficients in $i^* Rj_* (\bullet)$. Thanks to proper base change,
this definition is independent of the choice of $j$,
as is the long exact cohomology sequence 
\[
\quad\quad\quad\quad\quad\quad\quad\quad\quad \ldots \longto
H^n (X (\BC),\BQ) \otimes_\BQ F \longto 
\partial H^n (X (\BC),\BQ) \otimes_\BQ F \longto 
\]
\[
H^{n+1}_c (X (\BC),\BQ) \otimes_\BQ F \longto
H^{n+1} (X (\BC),\BQ) \otimes_\BQ F \longto 
\ldots \quad\quad\quad\quad\quad\quad\quad\quad
\]
(in the Hodge theoretic setting) resp. \
\[
\quad\quad\quad\quad\quad\quad\quad\quad \ldots \longto
H^n (X_{\bar{k}},\BQ_\ell) \otimes_\BQ F \longto
\partial H^n (X_{\bar{k}},\BQ_\ell) \otimes_\BQ F \longto
\]
\[
H^{n+1}_c (X_{\bar{k}},\BQ_\ell) \otimes_\BQ F \longto
H^{n+1} (X_{\bar{k}},\BQ_\ell) \otimes_\BQ F \longto
\ldots \quad\quad\quad\quad\quad\quad\quad\quad
\]
(in the $\ell$-adic setting). \\

Let us now assume in addition that an object $S \in Sm/k$ is given,
together with a factorization of the structure morphism of $X$, defining
on $X$ the structure of a proper, smooth scheme over $S$. 
In that situation, the \emph{Chow group} $\ch_{d_X} (X \times_S X)$
($d_X:=$ the absolute dimension of $X$) acts contravariantly
on the above exact cohomology sequences. By \cite[Thm.~2.2~(a)]{W6},
it also acts on the exact triangle $(*)$. 

\begin{Prop} \label{1Ea}
Let $X$ be proper and smooth over $S \in Sm/k$, and 
$e \in \ch_{d_X} (X \times_S X) \otimes_\BZ F$
an idempotent. Then
$H^n \! R \bigl( \dMgm(X)^e \bigr)$ is isomorphic to 
$\bigl( \partial H^n (X (\BC),\BQ) \otimes_\BQ F \bigr)^e$
(in the Hodge theoretic setting) resp. \
$\bigl( \partial H^n (X_{\bar{k}},\BQ_\ell) \otimes_\BQ F \bigr)^e$
(in the $\ell$-adic setting), for all $n$.
\end{Prop}

This result follows from \cite[Prop.~2.5]{W5}, if the class $e$ can
be represented by a cycle in $c_{1,2}(X,X) \otimes_\BZ F$.

\medskip

\begin{Proofof}{Proposition~\ref{1Ea}}
We imitate the proof of \cite[Prop.~2.5]{W5}, and
show that the image under $R$ of the
canonical morphism 
\[
\iota: \Mgm(X) \longto \Mcgm(X)
\]
can be $\ch_{d_X} (X \times_S X)$-equivariantly 
identified with the canonical morphism
\[
R \Gamma_c(X) \longto R \Gamma(X)
\]
in the target $D$ of $R$ of classes of complexes $R \Gamma_c(X)$ and 
$R \Gamma(X)$
computing cohomology with resp.\ without support. 

The first half of the proof of \cite[Prop.~2.5]{W5} 
contains the identification of $R(\iota)$
with $R \Gamma_c(X) \to R \Gamma(X)$. It thus remains to show that
this identification is compatible with the action of $\ch_{d_X} (X \times_S X)$.

Thus, let $[\FZ] \in \ch_{d_X} (X \times_S X)$ 
be the class of a cycle on $X \times_S X$. By \cite[Lemma~5.18]{Le},
it can be represented by a cycle $\FZ$ belonging to the group $c_S(X,X)$
of \emph{finite correspondences}, \emph{i.e.}, the
projection onto the first factor $X \times_S X \to X$ is finite on $\FZ$.
Similarly, the transposed class ${}^t [\FZ]$ can be represented by
$\FZ' \in c_S(X,X)$.
By \cite[Thm.~2.2~(a)]{W6}, the class $[\FZ]$ acts on $\Mgm(X)$ \emph{via}
the finite correspondence $\FZ$. By \cite[Thm.~2.2~(b), (c), Rem.~2.14]{W6},
it acts on $\Mcgm(X)$ \emph{via} the endomorphism dual to the one induced
by $\FZ'(d_X)[2d_X]$ under the duality isomorphism
\[
\Mgm (X)^*(d_X)[2d_X] \isoto \Mcgm (X)
\]
\cite[proof of Thm.~4.3.7~3]{V}. We thus identify the commutative diagram 
\[
\vcenter{\xymatrix@R-10pt{
        \Mgm (X) \ar[d]_{[\FZ]} \ar[r]^-{\iota} &
        \Mcgm (X) \ar[d]^{[\FZ]} \\
        \Mgm (X) \ar[r]^-{\iota} &
        \Mcgm (X)
\\}}
\]
with
\[
\vcenter{\xymatrix@R-10pt{
        \Mgm (X) \ar[d]_{\FZ} \ar[r]^-{\iota} &
        \Mgm (X)^*(d_X)[2d_X] \ar[d]^{(\FZ')^*(d_X)[2d_X]} \\
        \Mgm (X) \ar[r]^-{\iota} &
        \Mgm (X)^*(d_X)[2d_X]
\\}}
\]
By \cite[pp.~6--7]{DG}, $R$ sends $\FY$ to 
$\FY^* : R \Gamma (W) \to R \Gamma (V)$,
for any finite correspondence $\FY$ on the product $V \times_k W$
of two smooth $k$-schemes.   
It follows 
that $R$ sends the latter commutative diagram to the commutative diagram
\[
\vcenter{\xymatrix@R-10pt{
        R \Gamma (X) &
        R \Gamma (X)^*(-d_X)[-2d_X] \ar[l]  \\
        R \Gamma (X) \ar[u]^{\FZ^*} &
        R \Gamma (X)^*(-d_X)[-2d_X] \ar[l] \ar[u]_{((\FZ')^*)^*(-d_X)[-2d_X]}
\\}}
\] 
Now the endomorphism $((\FZ')^*)^*(-d_X)[-2d_X]$ of
$R \Gamma (X)^*(-d_X)[-2d_X]$ corres\-ponds to the endomorphism 
${}^t (\FZ')^* = {}^t [\FZ']^* = [\FZ]^*$ of
$R \Gamma_c (X)$. But this means precisely that our identification
of the image under $R$ of $\Mgm (X) \to \Mcgm (X)$ with 
the canonical morphism $R \Gamma_c (X) \to R \Gamma (X)$ is 
compatible with the action of $\ch_{d_X} (X \times_S X)$.
\end{Proofof}

Here is how the theory developed so far will be used in the sequel.

\begin{Cor} \label{1E}
Let $k$ be a field embeddable into $\BC$,
and $F$ a finite direct product of fields of characteristic zero. Let
$X$ be proper and smooth over $S \in Sm/k$, and 
$e \in \ch_{d_X} (X \times_S X) \otimes_\BZ F$
an idempotent. Assume that
the motive $\dMgm(X)^e$ lies in $\DFAbgM$,
and that 
$\bigl( \partial H^n (X (\BC),\BQ) \otimes_\BQ F \bigr)^e$
(in the Hodge theoretic setting) resp. \
$\bigl( \partial H^n (X_{\bar{k}},\BQ_\ell) \otimes_\BQ F \bigr)^e$
(in the $\ell$-adic setting)
is without weights $n - \beta,\ldots,n - (\alpha+1),n - \alpha$, 
for all $n \in \BZ$. Then  $\dMgm(X)^e$ 
is without weights $\alpha,\alpha+1,\ldots,\beta$.
If furthermore $\alpha \le -1$
and $\beta \ge 0$, then Assumption~\ref{1A} holds.
\end{Cor}

\begin{Proof}
Appply Theorem~\ref{1D} and Proposition~\ref{1Ea}.
\end{Proof}


\bigskip

%
%

\section{Statement of the main result}
\label{2}



In order to state our main result (Theorem~\ref{2Main}),
let us introduce the geometrical situation we are going to consider
from now on. $F$ is now a quadratic imaginary number field,
and the base $k$ equals $F$. The scheme
$X$ is a power of the universal Abelian threefold over a
Picard surface, and $e$ is associated  
to a \emph{dominant weight} $(k_1,k_2,c,d) \in \BZ^4$,
$k_1 \ge k_2 \ge 0$ 
(see below for the precise definition). 
Theorem~\ref{2Main} implies 
that in this context, Assumption~\ref{1A}
is satisfied if and only if $(k_1,k_2,c,d)$ is \emph{regular}: 
$k_1 > k_2 > 0$. More precisely,   
$\dMgm(X)$ lies in $\DFAbgM$, and setting $k := \min (k_1-k_2, k_2)$,
the $e$-part of $\dMgm(X)$
is without weights $-k,-(k-1),\ldots,k-1$. 
We then list the main consequences of 
this result (Corollaries~\ref{2G}--\ref{2I}), applying the theory developed in 
\cite[Sect.~4]{W4}.
The proof of Theorem~\ref{2Main} 
will be given in Section~\ref{3}. It is an application of
Corollary~\ref{1E}; in order to verify the hypotheses of the latter,
we heavily rely on the main result of \cite{Anc}. \\
 
Denote by $f \mapsto \bar{f}$ the non-trivial element of the Galois group
of $F$ over $\BQ$.
Note that there is an isomorphism of $\BQ$-algebras
\[
F \otimes_\BQ F \isoto F \times F \; , \; 
q \otimes f \longmapsto (q \cdot f, \bar{q} \cdot f) \; ;
\]
it is $F$-linear with respect to multiplication on the right on 
$F \otimes_\BQ F$. For any $F$-vector space $V'$,
this isomorphism induces a decomposition
\[
V' \otimes_\BQ F = V_+' \oplus V_-' 
\]
of the $F$-vector space $V' \otimes_\BQ F$
into two sub-spaces, where
\[
V_+' := \{ v \in V' \otimes_\BQ F \; , \; 
(f \otimes 1)v = (1 \otimes f) v \} \; ,
\]
\[
V_-' := \{ v \in V' \otimes_\BQ F \; , \; 
(f \otimes 1)v = (1 \otimes \bar{f}) v \} \; .
\] 
The projections $\pi_+$ and $\pi_-$ onto $V_+'$ and $V_-'$, respectively,
are induced by scalars in $F \otimes_\BQ F$: for any non-zero element $x \in F$
satisfying $\bar{x} = -x$, we have
\[
\pi_+ = \halb \bigl(1 \otimes 1 + x \otimes \frac{1}{x} \bigr) \; , \;
\pi_- = \halb \bigl(1 \otimes 1 - x \otimes \frac{1}{x} \bigr) \; .
\]
Conjugation $v \otimes f \mapsto v \otimes \bar{f}$ on $V' \otimes_\BQ F$
exchanges $V_+'$ and $V_-'$.
The restriction of $\pi_+$ to $V' = V' \otimes_\BQ \BQ \subset V' \otimes_\BQ F$
induces an isomorphism $V' \isoto V'_+$ of $\BQ$-vector spaces, which is
$F$-linear (for the given $F$-structure on $V'$). The analogous statement
holds for $\pi_-$, except that the induced isomorphism $V' \isoto V'_-$
is $F$-antilinear. Under these identifications, conjugation $V_+' \to V_-'$
corresponds to the identity on $V'$. \\  

Now fix a three-dimensional $F$-vector space $V$, together with an $F$-valued
Hermitian form $J$ on $V$ of signature $(2,1)$. 

\begin{Def}[{\cite[2.1]{G}}] \label{2A}
The group scheme $G$ over $\BQ$
is defined as the group of unitary similitudes
\[
G := GU(V,J) \subset \ReF \GL_F(V) \; .
\]
\end{Def}

Thus, for any $\BQ$-algebra $R$, the group $G(R)$ equals
\[
\{ g \in \GL_{F \otimes_\BQ R} (V \otimes_\BQ R) \; , \; 
\exists \, \lambda(g) \in R^* \; , \; 
J(g \argdot,g \argdot) = \lambda(g) \cdot J(\argdot,\argdot) \} \; .
\] 
In particular, the similitude norm $\lambda(g)$ defines a canonical morphism
\[
\lambda: G \longto \Gm \; .
\]
The decomposition $V \otimes_\BQ F = V_+ \oplus V_-$, 
together with the morphism $\lambda$, induces an isomorphism
\[
\phi: G_F := G \times_\BQ F \isoto \GL_F(V) \times_F \GFm \; , \;
g \longmapsto \bigl( g_{\tei V_+} , \lambda(g) \bigr)
\]
\cite[(2.1.1)]{G}; here we use the above identification
of $V$ and $V_+$. One deduces that $\lambda$ is an epimorphism
of algebraic groups. The reader will check that under $\phi$,
conjugation $f \mapsto \bar{f}$ on $G_F$ corresponds to the
involution $(h,t) \mapsto (\bar{t} h^*, \bar{t})$. Here, the map 
$h \mapsto h^*$ associates to $h \in \GL_F(V)$ the adjoint automorphism
$h^*$, characterized by the formula
\[
J(v_1,h^*v_2) = J(h^{-1}v_1,v_2) \; , \; \forall \, v_1, v_2 \in V \; .
\]
Under $\phi$, the group of $\BQ$-rational points of $G$ is identified
with
\[
\{ (h,q) \in \GL_F(V) \times \BQ^* \; , \; h = qh^* \} 
\]
\cite[(2.1.2)]{G}. From the definition, and from the isomorphism $\phi$,
one also deduces the following statement. 

\begin{Prop} \label{2B}
(a)~The group $G$ is split over $F$, but not over $\BQ$. \\[0.1cm]
(b)~The center $Z(G)$ of $G$ equals $\ReF \GFm \subset \ReF \GL_F(V)$
(inclusion of scalar automorphisms). In particular, it is isogeneous
to the product of $\Gm$ and a torus of compact type.
\end{Prop}

Thus, the irreducible algebraic representations of $G_F$ are indexed
by the dominant weights of $\GL_F(V) \times_F \GFm$, which we describe now.
$F$-split tori, together with an inclusion into a Borel subgroup of 
$\GL_F(V)$, are in bijection with $F$-bases of $V$; fix one such
basis, use it to identify $\GL_F(V)$ with $\GL_{3,F}$, the split torus
with the subgroup $T$ of diagonal matrices, and the Borel subgroup
with the subgroup of upper triangular matrices in $\GL_{3,F}$. 
The Weyl group is the symmetric group $\FS_3$; it acts by permuting
the base vectors, hence the coordinate functions of $T$.
Let us normalize things as follows: 
we consider quadruples $(k_1,k_2,c,r) \in \BZ^4$ satisfying 
two congruence relations: 
\[
c \equiv k_1 + k_2 \!\!\! \mod 2 \; ,\; 
r \equiv \frac{c+k_1+k_2}{2} \!\! \mod 2 \; .
\] 
To such a quadruple, let us associate the 
(representation-theoretic) weight 
\[
\alpha(k_1,k_2,c,r) : T \times_F \GFm \longto \GFm \; , 
\quad\quad\quad\quad\quad\quad\quad\quad
\]
\[
\quad\quad\quad\quad\quad\quad 
\bigl(  \diag (a,a^{-1}b,b^{-1}\nu) , f \bigr) \longmapsto
a^{k_1 - k_2} b^{k_2} \nu^{\frac{c-(k_1+k_2)}{2}} 
f^{-\halb (r + \frac{3c-(k_1+k_2)}{2})} \; .
\]
Note that
restriction of $\alpha(k_1,k_2,c,r)$ to $(T \cap \SL_{3,F}) \times_F \{1\}$
corresponds to the projection onto $(k_1,k_2)$. In particular,
the weight $\alpha(k_1,k_2,c,r)$ is dominant if and only if $k_1 \ge k_2 \ge 0$;
it is regular if and only if $k_1 > k_2 > 0$. Note also that the composition
of $\alpha(k_1,k_2,c,r)$ with the cocharacter
\[
\GFm \longto T \times_F \GFm \; , \; 
x \longmapsto \bigl(  \diag (x,x,x) , x^2 \bigr)
\]
equals 
\[
\GFm \longto \GFm \; , \; x \mapsto x^{-r} \; .
\]
The determinant on $T$ equals $\alpha(0,0,2,-3)$, and
$\lambda = \alpha(0,0,0,-2)$.
\forget{
\begin{Rem}
Conjugation on $G_F$ acts as follows on dominant weights:
\[
\alpha(k_1,k_2,c,r) \longmapsto \alpha(k_1,k_1-k_2,k_1-c,r) \; .
\]
\end{Rem}
}

\begin{Def} \label{2C}
The analytic space $\CH$ is defined as the complex open $2$-ball,
\[
\CH := \{
(z_1,z_2) \in \BC \; , \; \tei z_1 \tei^2 + \tei z_2 \tei^2 < 1 \} \; .
\]
\end{Def}

In order to identify  
the points of $\CH$ with certain morphisms 
of algebraic groups from the Deligne torus $\BS$ to $G_\BR$
\cite[Def.~2.5]{Anc}, one needs to choose an embedding of $F$ into $\BC$.
From now on, we thus assume such a choice of embedding to be fixed.
Therefore, we get a canonical action of $G(\BR)$ on $\CH$
(by conjugation). This action 
is by analytical automorphisms, and it is transitive \cite[(1.3.3)]{G}.
In fact, $(G,\CH)$ are \emph{pure Shimura data}
\cite[Def.~2.1]{P}. Their \emph{reflex field} \cite[Sect.~11.1]{P}
equals $F$ \cite[Lemma~4.2]{G}. According to Proposition~\ref{2B}~(b), 
the Shimura data $(G,\CH)$ satisfy condition 
$(+)$ from \cite[Sect.~5]{W3}. \\
  
Let us now fix additional data:
(A)~an open compact subgroup $K$ of $G(\BA_f)$
which is neat \cite[Sect.~0.6]{P}, 
(B)~a quadruple $(k_1,k_2,c,r) \in \BZ^4$ 
satisfying the above congruences
\[
c \equiv k_1 + k_2 \!\!\! \mod 2 \quad \text{and} \quad 
r \equiv \frac{c+k_1+k_2}{2} \!\! \mod 2 \; , 
\]   
and in addition,
\[
k_1 \ge k_2 \ge 0 \; .
\]
In other words, the character $\ua := \alpha(k_1,k_2,c,r)$ is dominant. 
We also impose the condition
\[
r \ge 0 \; .
\]
These data (A)--(B) are used as follows.
The \emph{Shimura variety}
$S := S^K (G,\CH)$ 
is an object of $Sm / F$. This is the \emph{Picard
surface} of level $K$ associated to $F$. 
According to Miyake (see \cite[Prop.~4.6, Thm.~4.9]{G}), it admits an interpretation as modular space of Abelian threefolds
with additional structures, among which a
complex multiplication by the ring of integers $\Fo_F$.
In particular, there is a universal family $A$ 
of Abelian threefolds over $S$. \\

Denote by $\Rep (G_F)$ the Tannakian category of algebraic
representations of $G_F$ in finite dimensional $F$-vector spaces.
The following result holds in the general context of Shimura varieties
of $PEL$-type.

\begin{Thm}[{\cite[Thm.~8.6]{Anc1}}] \label{2D}
Let $R \in \{ \BQ, F \}$.
There is an $R$-linear tensor functor
\[
\widetilde{\mu} = \widetilde{\mu}_R : \Rep (G_R) \longto CHM^s(S)_R \; .
\]
It has the following properties.
\begin{enumerate}
\item[(i)] The composition of $\widetilde{\mu}$
with the cohomological
Hodge theoretic reali\-zation is isomorphic to the \emph{canonical construction} functor $\mu$ 
(e.g.\ \cite[Thm.~2.2]{W1}) to the category
of admissible graded-polarizable variations of Hodge structure on $S$.
\item[(ii)] The functor $\widetilde{\mu}$ commutes with Tate twists.
\item[(iii)] The functor $\widetilde{\mu}$
maps the representation $V$
to the dual of the relative Chow motive $h^1(A/S)$.
\end{enumerate} 
\end{Thm}

Given that the representation on $V$ is faithful, it follows that
any relative Chow motive in the image of $\widetilde{\mu}$ 
is isomorphic to a direct sum of
direct factors of Tate twists of $h(A^{r_i}/S)$, for suitable $r_i \in \BN$,
where $A^{r_i}$ denotes the 
$r_i$-fold fibre product of $A$ over $S$. 
From the construction \cite[Thm.~8.6]{Anc1}, 
it follows that up to isomorphism, $\widetilde{\mu}_F$ is obtained from $\widetilde{\mu}_\BQ$
by formally tensoring with $F$ over $\BQ$, and passing to the pseudo-Abelian
completions. \\

\begin{Def} \label{2E}
(a)~Denote by $V_{\ua} \in \Rep (G_F)$ the irreducible representation
of highest weight $\ua$. \\[0.1cm]
(b)~Define ${}^{\ua} \CV \in \CHFSM$ as 
\[
{}^{\ua} \CV := \widetilde{\mu}_F(V_{\ua}) \; .
\]
\end{Def}

It will be useful to compare our parametrization of highest weights to that
of \cite[Sect.~4]{Anc}. There, the standard basis of characters of
the split torus $T \times_F \GFm$ is used. First step: write 
$(a,b,\gamma,d)$ instead of $(a,b,c,d)$ as in \loccit ,
and $\lambda(a,b,\gamma,d)$ for the associated character. Second step:
identify the change of parameters $(a,b,\gamma,d) \leftrightarrow
(k_1,k_2,c,r)$. 

\begin{Lem} \label{2F}
The character
$\alpha(k_1,k_2,c,r)$ of $T \times_F \GFm$ equals
\[ 
\lambda \bigl( \frac{c+k_1-k_2}{2},\frac{c-k_1+k_2}{2},\frac{c-(k_1+k_2)}{2},
-\halb (r + \frac{3c-(k_1+k_2)}{2}) \bigr) \; .
\]
The character $\lambda(a,b,\gamma,d)$ equals
\[
\alpha \bigl( a-\gamma,b-\gamma,a+b,-(a+b+\gamma+2d) \bigr) \; .
\]
\end{Lem}

\begin{Proof}
This follows directly from the definitions. We leave the computation
to the reader.
\end{Proof}

According to \cite[Lemma~4.13]{Anc}, the Chow motive
${}^{\ua} \CV$ is a direct factor of $h^r(A^N/S)$,
for a suitable large enough positive integer $N$. 
More precisely, the representation $V_{\ua}$ is a direct factor of the 
exterior algebra $\Lambda^\bullet (V^N)^\vee$ 
of the dual of the representation $V^N$. Therefore, there is an idempotent
$e_{\ua}$ acting on $\Lambda^\bullet (V^N)^\vee$, whose image equals $V_{\ua}$.
Applying $\widetilde{\mu}_F$, we get an idempotent
acting on the relative Chow motive $h(A^N/S)$, 
equally denoted by $e_{\ua}$,
and whose image equals ${}^{\ua} \CV$. According to 
\cite[Thm.~2.2~(a)]{W6}, the relative Chow motive
${}^{\ua} \CV$
gives rise to an exact triangle
\[
\dMgm \bigl({}^{\ua} \CV \bigr) \longto \Mgm \bigl({}^{\ua} \CV \bigr) 
\stackrel{u}{\longto} \Mcgm \bigl({}^{\ua} \CV \bigr) 
\longto \dMgm \bigl({}^{\ua} \CV \bigr)[1] 
\]
in $\DFFgM$. By functoriality, and by \cite[Thm.~2.2~(a1)]{W6},
this triangle coincides with the $e_{\ua}$-part of the triangle
\[
\dMgm(A^N) \longto \Mgm(A^N) 
\longto \Mcgm(A^N) \longto \dMgm(A^N)[1]
\]
denoted $(*)$ in Section~\ref{1}.

\begin{Rem}
The restriction ``$r \ge 0$'' on our character $\ua = \alpha(k_1,k_2,c,r)$
is not very serious. For negative $r$, the Chow motive ${}^{\ua} \CV$
is defined the same way, the only difference being that it is not effective.
Indeed, it is then a direct factor of  
$h^{-r}(A^N/S)(-r)$, for some $N >> 0$.
\end{Rem}

Here is our main result.

\begin{Thm} \label{2Main}
The boundary motive $\dMgm(A^N)$ lies in the
triangulated sub-category $\DFQAbgM$ of $\DFQgM$. Its direct factor
$\dMgm({}^{\ua} \CV)$ is without weights 
\[
-k, -(k-1), \ldots, k-1 \; ,
\]
where $k := \min (k_1-k_2, k_2)$. Both weights
$-(k+1)$ and $k$ \emph{do} occur in $\dMgm({}^{\ua} \CV)$.
In particular, Assumption~\ref{1A} holds 
for $\dMgm({}^{\ua} \CV)$ if and only if $k \ge 1$, \emph{i.e.}, 
if and only if $\ua$ is regular.
\end{Thm}

The category $\DFFAbgM$ being pseudo-Abelian (Corollary~\ref{1Ja}),
the motive $\dMgm({}^{\ua} \CV)$ belongs to $\DFFAbgM$.
Theorem~\ref{2Main} should be compared to \cite[Thm.~3.5]{W5},
which treats the case of Hilbert-Blumenthal varieties. 
As here, regularity of the character is sufficient for the
corresponding direct factor of the boundary motive to avoid
weights $-1$ and $0$. However, as soon as the dimension of
the Hilbert--Blumenthal variety is strictly greater than one,
there are irregular dominant characters whose associated 
direct factors nonetheless avoid these weights. \\

Theorem~\ref{2Main}
will be proved in Section~\ref{3}.  
Let us give its main corollaries,
assuming that $k = \min (k_1-k_2, k_2) \ge 1$, \emph{i.e.}, $k_1 > k_2 > 0$.
Consider the weight filtration 
\[ 
C_{\le -(k+1)} \longto \dMgm \bigl( {}^{\ua} \CV \bigr) \longto 
C_{\ge k} \longto C_{\le -(k+1)}[1]
\]
avoiding weights $-k, \ldots, k-1$ \cite[Def.~1.6, Cor.~1.9]{W4}. 
Thus, 
\[
C_{\le -(k+1)} \in DM_{gm}^{Ab}(F)_{F,w \le -(k+1)} \; , \;
C_{\ge k} \in DM_{gm}^{Ab}(F)_{F,w \ge k} \; .
\]
Furthermore, according to Theorem~\ref{2Main},
\[
C_{\le -(k+1)} \not \in DM_{gm}^{Ab}(F)_{F,w \le -(k+2)} \; , \;
C_{\ge k} \not \in DM_{gm}^{Ab}(F)_{F,w \ge k+1} \; .
\] 

\begin{Cor}[{\cite[Thm.~4.3]{W4}}]  \label{2G}
Assume $k_1 > k_2 > 0$, \emph{i.e.}, $k \ge 1$. \\[0.1cm]
(a)~The motive $\Mgm({}^{\ua} \CV)$ is without weights $-k,-(k-1), \ldots, -1$,
and the motive $\Mcgm({}^{\ua} \CV)$ is without weights $1,2, \ldots, k$. 
The Chow motives $\Gr_0 \Mgm({}^{\ua} \CV)$ 
and $\Gr_0 \Mcgm({}^{\ua} \CV)$ \cite[Prop.~2.2]{W4} are defined,
and they behave functorially with respect to $\Mgm({}^{\ua} \CV)$ and 
$\Mcgm({}^{\ua} \CV)$. In particular, any endomorphism of 
$\Mgm({}^{\ua} \CV)$ induces an endomorphism of $\Gr_0 \Mgm({}^{\ua} \CV)$,
and any endomorphism of $\Mcgm({}^{\ua} \CV)$ induces an endomorphism
of $\Gr_0 \Mcgm({}^{\ua} \CV)$. \\[0.1cm]
(b)~There are canonical exact triangles 
\[
C_{\le -(k+1)} \longto \Mgm({}^{\ua} \CV) \stackrel{\pi_0}{\longto}
\Gr_0 \Mgm({}^{\ua} \CV) \longto C_{\le -(k+1)}[1]
\]
and
\[
C_{\ge k} \longto \Gr_0 \Mcgm({}^{\ua} \CV) \stackrel{i_0}{\longto}
\Mcgm({}^{\ua} \CV) \longto C_{\ge k}[1] \; ,
\]
which are stable under the natural action
of the endomorphism rings of $\Mgm({}^{\ua} \CV)$ and $\Mcgm({}^{\ua} \CV)$, 
respectively. In particular, weight $-(k+1)$ occurs in $\Mgm({}^{\ua} \CV)$,
and weight $k+1$ occurs in $\Mcgm({}^{\ua} \CV)$. \\[0.1cm]
(c)~There is a canonical isomorphism 
$\Gr_0 \Mgm({}^{\ua} \CV) \isoto \Gr_0 \Mcgm({}^{\ua} \CV)$ in 
$CHM(F)_F$. 
As a morphism, it is uniquely determined by the property
of making the diagram
\[
\vcenter{\xymatrix@R-10pt{
        \Mgm({}^{\ua} \CV) \ar[r]^-u \ar[d]_{\pi_0} &
        \Mcgm({}^{\ua} \CV) \\
        \Gr_0 \Mgm({}^{\ua} \CV) \ar[r] &
        \Gr_0 \Mcgm({}^{\ua} \CV) \ar[u]_{i_0}
\\}}
\]
commute. \\[0.1cm]
(d)~Let $N \in CHM(F)_F$ be a Chow motive. Then $\pi_0$ and $i_0$ induce
isomorphisms
\[
\Hom_{CHM(F)_F} \bigl( \Gr_0 \Mgm({}^{\ua} \CV) , N \bigr) \isoto
\Hom_{\DFFgM} \bigl( \Mgm({}^{\ua} \CV) , N \bigr)
\]
and
\[
\Hom_{CHM(F)_F} \bigl( N , \Gr_0 \Mcgm({}^{\ua} \CV) \bigr) \isoto
\Hom_{\DFFgM} \bigl( N , \Mcgm({}^{\ua} \CV) \bigr) \; .
\]
(e)~Let $\Mgm({}^{\ua} \CV) \to N \to \Mcgm({}^{\ua} \CV)$ 
be a factorization of the canonical morphism 
$u: \Mgm({}^{\ua} \CV) \to \Mcgm({}^{\ua} \CV)$ through a Chow motive
$N \in CHM(F)_F$. 
Then $\Gr_0 \Mgm({}^{\ua} \CV) = \Gr_0 \Mcgm({}^{\ua} \CV)$ is canonically
a direct factor of $N$, with a canonical direct complement.
\end{Cor} 

\begin{Proof}
This is \cite[Thm.~4.3]{W4}. The functoriality statements of
parts (a) and (b) are more general than \loccit , and follow from  
\cite[Prop.~1.7]{W4}.
\end{Proof}

Henceforth, we identify $\Gr_0 \Mgm({}^{\ua} \CV)$ and 
$\Gr_0 \Mcgm({}^{\ua} \CV)$
\emph{via} the canonical isomorphism of Corollary~\ref{2G}~(c).
The equivariance statements from
Corollary~\ref{2G}~(a), (b) apply in particular
to cycles coming from the \emph{Hecke algebra} associated to
the Shimura variety $S$.

\begin{Cor} \label{2H}
Assume $k \ge 1$. Then $\Gr_0 \Mgm({}^{\ua} \CV)$
carries a natural action of the Hecke algebra $\FH(K,G(\BA_f))$
associated to the neat open compact subgroup $K$ of $G(\BA_f)$.
\end{Cor}

\begin{Proof}
This is an application of \cite[Ex.~2.16]{W6}.
We refer to the proof of \cite[Cor.~3.8]{W5},
where the details of the construction are spelled out
for Hilbert--Blumenthal varieties. The reader
will have no difficulties to translate them to the present context. 
\end{Proof}

\begin{Cor}[{\cite[Cor.~4.6]{W4}}]  \label{2I} 
Assume $k \ge 1$, and
let $\widetilde{A^N}$ be any smooth compactification of $A^N$. Then
$\Gr_0 \Mgm({}^{\ua} \CV)$ is canoni\-cally
a direct factor of the Chow motive
$\Mgm(\widetilde{A^N})$, with a canonical direct complement.
\end{Cor}

Furthermore, \cite[Thm.~4.7, Thm.~4.8]{W4} on the Hodge theoretic
and $\ell$-adic realizations \cite[Cor.~2.3.5, Cor.~2.3.4 and Corrigendum]{Hu}
apply, and tell us in particular that $\Gr_0 \Mgm({}^{\ua} \CV)$
is mapped to the part of \emph{interior cohomology} of $A^N$
fixed by $e_{\ua}$. In particular, the $L$-function of the 
Chow motive $\Gr_0 \Mgm({}^{\ua} \CV)$
is computed \emph{via} (the $e_{\ua}$-part of) interior cohomology of $A^N$.

\begin{Def}[{\cite[Def.~4.9]{W4}}] \label{2J}
Let $k \ge 1$. 
We call $\Gr_0 \Mgm({}^{\ua} \CV)$ the \emph{$e_{\ua}$-part of the
interior motive of $A^N$}.
\end{Def}

\begin{Rem} \label{2K}
By \cite[Thm.~4.14]{W4}, control of the reduction of \emph{some}
compactification of $A^N$ implies control of certain properties
of the $\ell$-adic realization of $\Gr_0 \Mgm({}^{\ua} \CV)$.
To the best of the author's knowledge, the sharpest result known 
about reduction of compactifications of $A^N$ is stated
in \cite[Sect.~1.2]{R}: there exists a smooth compactification
of $A^N$ having good reduction at each prime ideal $\Fp$
dividing neither the level $n$ of $K$ nor the absolute discriminant
$d$ of $F$. Since the argument is somewhat involved, 
let us reproduce it: first, according to 
\cite[Theorems on p.~34, Sect.~8]{La}, the Picard surface $S$
and its canonical smooth toroidal compactification $\bar{S}$
both admit models $\CS$ and $\bar{\CS}$, respectively, which are smooth
over the ring $\Fo_F[1/(nd)]$; furthermore, the $\Fo_F[1/(nd)]$-scheme
$\bar{\CS}$ is proper. In addition, the complement $\CD$ of $\CS$ in
$\bar{\CS}$ continues to be a smooth divisor, and
the universal Abelian scheme $A$ over $S$ admits canonical extensions,
first to an Abelian scheme $\CA$ over $\CS$, and then,
to a semi-Abelian scheme $\CG$ over $\bar{\CS}$. The reader may find it
useful to consult \cite[Chap.~I]{Bel}, in particular 
\cite[Thm.~I.3.2.4, Thm.~I.5.1.1, Prop.~I.5.3.4]{Bel} for the proofs
of these statements. Second, one applies \cite[Thm.~1]{R1}:
for any $N \ge 0$,
there exists a compactification $\bar{\CG^N}$ of $\CG^N$, which
is regular, and proper and flat over $\bar{\CS}$.
In fact, looking at the construction of \loccit, which uses
toric compactifications, one sees that the morphism
$\bar{\CG^N} \to \bar{\CS}$ is in fact semi-stable.
Writing \'etale-local equations (and using that $\CD$ is smooth
over $\Fo_F[1/(nd)]$!), one sees that $\bar{\CG^N}$ 
is smooth over $\Fo_F[1/(nd)]$.

\cite[Thm.~4.14]{W4} then yields the following: 
(a)~for all prime ideals $0 \ne \Fp$ not dividing $nd$, the $\Fp$-adic
realization of $\Gr_0 \Mgm({}^{\ua} \CV)$ is crystalline,
(b)~if furthermore $\Fp$ and $\ell$ are coprime, then
the $\ell$-adic realization of $\Gr_0 \Mgm({}^{\ua} \CV)$
is unramified at $\Fp$. 
\end{Rem}

\begin{Cor} \label{2L}
Let $0 \ne \Fp$ be a prime ideal of $\Fo_F$ dividing neither the level of $K$
nor the discriminant of $F$. 
Let $\ell$ be coprime to $\Fp$.
Then the characteristic polynomials of the following coincide: (1)~the action
of Frobenius $\phi$ on the $\phi$-filtered module associated to
the (crystalline) $\Fp$-adic realization of 
$\Gr_0 \Mgm({}^{\ua} \CV)$, (2)~the action of a geometric 
Frobenius automorphism at
$\Fp$ on the (unramified) $\ell$-adic realization of 
$\Gr_0 \Mgm({}^{\ua} \CV)$.  
\end{Cor}

\begin{Proof}
According to our construction, and what was recalled in Remark~\ref{2K},
there is a smooth and proper scheme $\bar{\CG^N_{\BF_\Fp}}$ 
over the residue field $\BF_\Fp$ of $\Fp$, and an endomorphism 
$e_{\ua,\Fp}$ of the Chow motive associated to $\bar{\CG^N_{\BF_\Fp}}$,
whose images in the endomorphism rings of the realizations of
$\bar{\CG^N_{\BF_\Fp}}$ become idempotent; furthermore, 
the latter are projectors onto
the realizations of $\Gr_0 \Mgm({}^{\ua} \CV)$. The claim thus follows from
\cite[Thm.~2.~2)]{KM}.
\end{Proof}

\begin{Rem}
The case $\ua = \alpha(0,0,0,0)$ is not covered by Corolla\-ries~\ref{2G}--\ref{2I}. 
It concerns
the motive and the motive with compact support of the base scheme $S$.
In this situation, the best replacement 
of the interior motive is the \emph{intersection motive}
$M^{!*}(S)$ of $S$ with respect to its Baily--Borel compactification 
\cite{CM,W7}.
\end{Rem}


\bigskip

%
%

\section{Proof of the main result}
\label{3}



We keep the notation of the preceding section. In order to prove
Theorem~\ref{2Main}, the idea is to apply the
criterion from Corollary~\ref{1E}. 

\begin{Thm} \label{3A}
For any integer $N \ge 0$,
the boundary motive $\dMgm(A^N)$ lies in the
triangulated sub-category $\DFQAbgM$ of $\DFQgM$.
\end{Thm}

\begin{Proof}
The variety $A^N$ is a \emph{mixed Shimura variety} over 
$S = S^K (G,\CH)$. More precisely, using the description
of the morphisms $h : \BS \to G_\BR$ from \cite[Def.~2.5]{Anc}, the
representation $V$ of $G$ is seen to be of Hodge type
$\{ (-1,0) , (0,-1) \}$ in the sense of \cite[Sect.~2.16]{P}.
The same statement is then true for the $r$-th power $V^N$
of $V$. By \cite[Prop.~2.17]{P}, this allows for the construction
of the \emph{unipotent extension} $(P^N,\FX^N)$ of $(G,\CH)$
by $V^N$. 

The pair $(P^N,\FX^N)$ constitute \emph{mixed
Shimura data} \cite[Def.~2.1]{P}. By construction, they come endowed with
a morphism $\pi^N: (P^N,\FX^N) \to (G,\CH)$ of Shimura data,
identifying $(G,\CH)$ with the pure Shimura data
underlying $(P^N,\FX^N)$. In particular, $(P^N,\FX^N)$ also satisfy
condition $(+)$ from \cite[Sect.~5]{W3}.

Now there is an open compact neat subgroup $K^N$ of $P^N(\BA_f)$,
whose image under $\pi^N$ equals $K$, and
such that $A^N$ is identified with the \emph{mixed Shimura variety}
$S^{K^N} (P^N,\FX^N)$ \cite[Sect.~3.22, Thm.~11.18 and 11.16]{P}.
Furthermore, the morphism $\pi^N$ of Shimura data induces a morphism
$S^{K^N} (P^N,\FX^N) \to S^K (G,\CH)$, which is identified with the
structure morphism of $A^N$.

In order to obtain control on
the boundary motive of $A^N$, we choose a smooth
\emph{toroidal compactification} $\widetilde{A^N}$. It is associated to a 
\emph{$K^N$-admissible complete smooth cone 
decomposition} $\FS$, \emph{i.e.}, a collection of subsets of
\[
\CC (P^N,\FX^N) \times P^N (\BA_f)
\]
satisfying the axioms of \cite[Sect.~6.4]{P}. Here,
$\CC (P^N,\FX^N)$
denotes the \emph{co\-ni\-cal complex} associated to $(P^N,\FX^N)$
\cite[Sect.~4.24]{P}. 

We refer to \cite[9.27, 9.28]{P} for criteria sufficient to guarantee
the existence of the associated compactification
$\widetilde{A^N} := S^{K^N} (P^N , \FX^N , \FS)$. 
It comes equipped with a natural (finite) 
stratification into locally closed strata.
The unique open stratum is $A^N$. Any stratum 
$\widetilde{A^N_\sigma}$ different from the generic
one is associated to a \emph{rational boundary component} $(P_1,\FX_1)$
of $(P^N,\FX^N)$ \cite[Sect.~4.11]{P} which is \emph{proper},
\emph{i.e.}, unequal to $(P^N,\FX^N)$. 

First, co-localization for the boundary motive \cite[Cor.~3.5]{W2}
tells us that $\dMgm(A^N)$ is a successive extension of 
(shifts of) objects of the form
\[
\Mgm (\widetilde{A^N_\sigma}, i_\sigma^! \, j_! \, \BZ) \; .
\]
Here, $j$ denotes the open immersion of $A^N$ into $\widetilde{A^N}$,
$i_\sigma$ runs through the immersions of 
the strata $\widetilde{A^N_\sigma}$ different from $A^N$ into $\widetilde{A^N}$,
and $\Mgm (\widetilde{A^N_\sigma}, i_\sigma^! \, j_! \, \BZ)$
is the motive of $\widetilde{A^N_\sigma}$ 
with coefficients in $i_\sigma^! \, j_! \, \BZ$
defined in \cite[Def.~3.1]{W2}. 

Next, by \cite[Thm.~6.1]{W3}, there is an isomorphism
\[
\Mgm (\widetilde{A^N_\sigma}, i_\sigma^! \, j_! \, \BZ) \isoto
\Hom (\BZ (\sigma),\Mgm (S^{K_1}(P_1,\FX_1))) [\dim \sigma]  \; .
\]
Recall \cite[p.~971]{W3} that the \emph{group of orientations}
$\BZ (\sigma)$ is (non-canonically) isomorphic to $\BZ$,
hence 
\[
\Hom (\BZ (\sigma),\Mgm (S^{K_1}(P_1,\FX_1))) \cong 
\Mgm (S^{K_1}(P_1,\FX_1))  \; .
\]
$S^{K_1}(P_1,\FX_1)$ is a Shimura variety associated to the data
$(P_1,\FX_1)$ and an open compact neat subgroup $K_1$ of $P_1(\BA_f)$. 
In order to show our claim, we are thus reduced to showing that
$\Mgm (S^{K_1})$ is an object of $\DFQAbgM$, 
for any Shimura variety
$S^{K_1} = S^{K_1}(P_1,\FX_1)$ associated to a proper rational
boundary component $(P_1,\FX_1)$ of $(P^N,\FX^N)$, and any
open compact neat subgroup $K_1$ of $P_1(\BA_f)$.

Given that $P^N$ is a unipotent extension of $G$,
the pure Shimura data underlying $(P_1,\FX_1)$ coincides
with the pure Shimura data underlying some proper rational
boundary component of $(G,\CH)$. 
These boundary components are determined in \cite[Sect.~3]{Anc}.
In particular \cite[Lemma~3.8, Prop.~3.12]{Anc}, the pure Shimura data
$(G_1,\CH_1)$ underlying any such component all equal
$(\ReF \GFm,\{*\})$, where $\{*\}$
is a single point. 
Altogether, we see that $(P_1,\FX_1)$
is a unipotent extension of $(\ReF \GFm,\{*\})$.

We are ready to conclude.
As follows directly
from the definition of the canonical model
(cmp.~\cite[Sect.~11.3, 11.4]{P}),
the pure Shimura variety $S^{\pi^N(K_1)}(\ReF \GFm,\{*\})$ 
underlying $S^{K_1}$ equals the spectrum of a finite (Abelian) field 
extension $C$ of $F$. 
Consider the factorization of $\pi^N: (P_1,\FX_1) \to (G_1,\CH_1)$
corresponding to the weight filtration $1 \subset U \subset W$ 
of the unipotent radical 
$W$ of $P_1$. It gives the following:
\[
\vcenter{\xymatrix@R-10pt{
(P_1,\FX_1) \ar@{->>}[r]^-{\pi_t} \ar@/_2pc/[rr]^-{\pi^N}&
(P'_1,\FX'_1) := (P_1,\FX_1)/U \ar@{->>}[r]^-{\pi_a} &
(G_1,\CH_1)
\\}}
\]
On the level of Shimura varieties, we get:
\[
\vcenter{\xymatrix@R-10pt{
S^{K_1} \ar@{->>}[r]^-{\pi_t} \ar@/_2pc/[rr]^-{\pi^N} &
S^{\pi_t (K_1)} = S^{\pi_t (K_1)}(P'_1,\FX'_1) \ar@{->>}[r]^-{\pi_a} &
S^{\pi^N(K_1)}(\ReF \GFm,\{*\})
\\}}
\]
By \cite[3.12--3.22~(a), Prop.~11.10]{P}, 
$\pi_a$ is in a natural way a torsor under an
Abelian variety, while $\pi_t$ is a torsor under a split torus.
After base change over $S^{\pi^N(K_1)} = S^{\pi^N(K_1)}(\ReF \GFm,\{*\})$ 
to the algebraic closure of $C$, the morphism $\pi_a$ thus
becomes isomorphic to an Abelian variety; in particular, 
\[
\Mgm \bigl( S^{\pi_t (K_1)} \bigr) \in \CHFQAbM \; .  
\]
As for $S^{K_1}$, it is the fibre product over $S^{\pi_t (K_1)}$ of
a finite number $m$ of $S^{\pi_t (K_1)}$-schemes, all of which are
isomorphic to complements $\FL_i^*$ of the zero section in line bundles
$\FL_i$ over $S^{\pi_t (K_1)}$, $i = 1,\ldots,m$. 

Let $Y$ be an $S^{\pi_t (K_1)}$-scheme, and $i \in \{ 1,\ldots,m \}$.
\emph{Homotopy invariance} and the \emph{Mayer-Vietoris axiom} \cite[Sect.~2.2]{V}
together imply that the canonical morphism 
\[
\Mgm \bigl( \FL_i \times_{S^{\pi_t (K_1)}} Y \bigr) \longto \Mgm (Y) 
\]
is an isomorphism. Furthermore, the \emph{Gysin triangle} 
\[
\Mgm (Y)(1)[1] \to
\Mgm \bigl( \FL_i^* \times_{S^{\pi_t (K_1)}} Y \bigr) \to  
\Mgm \bigl( \FL_i \times_{S^{\pi_t (K_1)}} Y \bigr) \to
\Mgm (Y)(1)[2]
\]
associated to the zero section of $\FL_i$
is exact \cite[Prop.~3.5.4]{V}.
Together, the two statements imply that if $\Mgm (Y)$ belongs to $\DFQAbgM$,
then so does $\Mgm (\FL_i^* \times_{S^{\pi_t (K_1)}} Y)$.
\end{Proof}

\begin{Rem} \label{3B}
A more detailed analysis of
the proof reveals that the boundary motive $\dMgm(A^N)$ actually lies in the
full triangulated sub-category of $\DFQAbgM$
generated by Tate twists and Chow motives $\Mgm(X)$
associated to (proper) schemes $X \in Sm/F$, whose base change
to $\bar{F}$
is isomorphic to a finite disjoint union of products of elliptic curves
with complex multiplication by $F$.
\end{Rem}

In order to apply the results from Section~\ref{1}, we need to identify
the weights occurring in the Hodge structure on the $e_{\ua}$-part 
of the boundary cohomology of $A^N$,
\[
\bigl( \partial H^n \bigl( A^N (\BC),\BQ \bigr) 
                   \otimes_\BQ F \bigr)^{e_{\ua}} \; ,
\]
for all integers $n$. 

\begin{Prop} \label{3C}
There is a canonical isomorphism of Hodge structures
\[
\bigl( \partial H^n \bigl( A^N (\BC),\BQ \bigr) \otimes_\BQ F \bigr)^{e_{\ua}} 
\isoto
\partial H^{n-r} \bigl( S (\BC), \mu(V_{\ua}) \bigr) 
\]
for all integers $n$.
\end{Prop}

\begin{Proof}
Recall that
the representation $V_{\ua}$ is a direct factor,
\emph{via} the idempotent $e_{\ua}$, of the 
exterior algebra $\Lambda^\bullet (V^N)^\vee$ 
of the dual of the representation $V^N$. 
Applying the cohomological Hodge theoretic realization
$\mu$, we get an idempotent
acting on all 
relative cohomology objects $\CH^q(A^N/S)$,
whose image equals $\mu(V_{\ua})$ if $q=r$, and zero otherwise.

Consider the spectral sequence of Hodge structures
\[
E_2^{p,q} = \partial H^p \bigl( S (\BC), \CH^q(A^N/S) \bigr) \Longrightarrow 
\partial H^{p+q} \bigl( A^N (\BC),\BQ \bigr) \; .
\]
It is compatible with the action of the Chow group
\[
\ch_{3N+2} (A^N \times_S A^N) \otimes_\BZ F = \End_{CHM^s(S)_F}(h(A^N/S)) \; .
\]
Note that on the $E_2$-terms of the spectral sequence, the Chow group
acts only on $\CH^\bullet(A^N/S)$. Thus, for the fixed part under $e_{\ua}$, 
we obtain a spectral sequence
\[
\partial H^p \bigl( S (\BC), \CH^q(A^N/S)^{e_{\ua}} \bigr) \Longrightarrow 
\bigl( \partial H^{p+q} 
              \bigl( A^N (\BC),\BQ \bigr) \otimes_\BQ F \bigr)^{e_{\ua}} \; .
\]
But $\CH^q(A^N/S)^{e_{\ua}} = 0$ for $q \ne r$, while
$\CH^r(A^N/S)^{e_{\ua}} = \mu(V_{\ua})$.
\end{Proof}

\begin{Rem} \label{3D}
The same proof show that 
\[
\bigl( H^n \bigl( A^N (\BC),\BQ \bigr) \otimes_\BQ F \bigr)^{e_{\ua}} 
\isoto
H^{n-r} \bigl( S (\BC), \mu(V_{\ua}) \bigr) 
\]
for all integers $n$, and likewise for cohomology wih compact support. 
Saper's vanishing theorem 
\cite[Thm.~5]{Sa} says that if $\ua$ is regular, then the groups
$H^{n-r} \bigl( S (\BC), \mu(V_{\ua}) \bigr)$ vanish for $n < r+2$.
By duality, one abtains that 
$H^{n-r}_c \bigl( S (\BC), \mu(V_{\ua}) \bigr) = 0$ for $n > r+2$.
It follows that interior cohomology
\[
\bigl( H^n_! \bigl( A^N (\BC),\BQ \bigr) \otimes_\BQ F \bigr)^{e_{\ua}} 
\isoto
H^{n-r}_! \bigl( S (\BC), \mu(V_{\ua}) \bigr) 
\]
vanishes unless $n = r+2$. By comparison, the analogous statement
is true for $\ell$-adic cohomology. In other words, 
the realizations of $\Gr_0 \Mgm({}^{\ua} \CV)$
are concentrated in the single cohomological degree $r+2$
if $\ua$ is regular, and they take the values
\[
\bigl( H^{r+2}_! \bigl( A^N (\BC),\BQ \bigr) \otimes_\BQ F \bigr)^{e_{\ua}} 
\isoto
H^2_! \bigl( S (\BC), \mu(V_{\ua}) \bigr) 
\]
(in the Hodge theoretic setting) resp.\ 
\[
\bigl( H^{r+2}_! \bigl( A^N_{\bar{F}},\BQ_\ell \bigr) 
\otimes_\BQ F \bigr)^{e_{\ua}} \isoto
H^2_! \bigl( S_{\bar{F}}, \mu_\ell(V_{\ua}) \bigr) 
\]
(in the $\ell$-adic setting); see e.g.\ \cite[Sect.~(4.1)]{P2}
for the $\ell$-adic version $\mu_\ell$ of the canonical construction functor. 
\end{Rem}

\begin{Proofof}{Theorem~\ref{2Main}}
According to Proposition~\ref{3C}, 
we need to control $\partial H^{n-r} ( S (\BC), \mu(V_{\ua}) )$. 
We use the \emph{Baily--Borel compactification} 
$S^*$ of $S$. The complement of $S$ consists
of finitely many cusps;
the boundary cohomology of $S (\BC)$ therefore coincides with the direct sum
over the cusps of the \emph{degeneration} of the coefficients to the 
boundary of $S^*$. 

The first claim of Theorem~\ref{2Main}
is Theorem~\ref{3A}.
Fix a cusp of $S^*$, 
denote by
$j$ the open immersion of $S$,
and by $i$ the closed immersion of the cusp into $S^*$. 
We need to know the weights occurring in 
\[
i^* R^{n-r} j_* \bigl(  \mu(V_{\ua}) \bigr) \; .
\]
Recall (Lemma~\ref{2F}) that 
\[ 
\ua =
\lambda \bigl( \frac{c+k_1-k_2}{2},\frac{c-k_1+k_2}{2},\frac{c-(k_1+k_2)}{2},
-\halb (r + \frac{3c-(k_1+k_2)}{2}) \bigr) 
\]
in the parametrization of \cite[Sect.~4]{Anc}. The main result of \cite{Anc}
implies that 
\begin{enumerate}
\item[(o)] $0 \ne i^* R^0 j_* (\mu(V_{\ua}))$ is of weight $r - k_1$,
\item[(i)] $0 \ne i^* R^1 j_* (\mu(V_{\ua}))$ is of weights $(r+1) - k_2$ and
$(r+1) - (k_1-k_2)$,
\item[(ii)] $0 \ne i^* R^2 j_* (\mu(V_{\ua}))$ is of weights 
$(r+2) + k_2 + 1$ and $(r+2) + (k_1 - k_2) + 1$,
\item[(iii)] $0 \ne i^* R^3 j_* (\mu(V_{\ua}))$ is of weight 
$(r+3) + k_1 + 1$,
\end{enumerate} 
and that $i^* R^m j_* (\mu(V_{\ua})) = 0$ whenever $m < 0$ or $m > 3$
\cite[Thm.~1.2]{Anc}.

Therefore, $\partial H^{n-r} ( S (\BC), \mu(V_{\ua}) )$ is without weights
$n-(k-1),\ldots,n+k-1,n+k$, where $k = \min (k_1-k_2, k_2)$,
while weight $n-k$ occurs for $n = r+1$, and weight $n+k+1$
for $n = r+2$.
Our claim thus follows from Proposition~\ref{3C} and Corollary~\ref{1E}.
\end{Proofof}


\bigskip

%
%

\section{The motive for an automorphic form}
\label{4}



This final section contains the analogues for Picard surfaces
of the main results from \cite{Sc}. Since we shall not restrict
ourselves to the case of Hecke eigenforms, our
notation becomes a little more technical than in \loccit ; 
we chose however to keep it compatible with the one used in \cite{Editors}. \\

We continue to consider the situation 
of Sections~\ref{2} and \ref{3}. In particular,
we fix a dominant $\ua = \alpha(k_1,k_2,c,r)$, which we
assume to be regular, \emph{i.e.}, $k_1 > k_2 > 0$. Consider 
the Chow motive $\Gr_0 \Mgm({}^{\ua} \CV)$. According to Remark~\ref{3D},
its Hodge theoretic realization equals
\[
R \bigl( \Gr_0 \Mgm({}^{\ua} \CV) \bigr) = 
H^2_! \bigl( S (\BC), \mu(V_{\ua}) \bigr)[-(r+2)]
\; .
\]
By Corollary~\ref{2H}, the Hecke algebra $\FH(K,G(\BA_f))$ acts on
$\Gr_0 \Mgm({}^{\ua} \CV)$. 

\begin{Thm}[{\cite[Chap.~2, Thm.~2, p.~50]{Ha}}]  \label{4A}
Let $L$ be any field extension of $F$. Then the 
$\FH(K,G(\BA_f)) \otimes_F L$-module 
$H^2_! ( S (\BC), \mu(V_{\ua}) ) \otimes_F L$
is semi-simple. 
\end{Thm}

Note that \cite[Chap.~3, Sect.~4.3.5]{Ha} gives a proof of
Theorem~\ref{4A}, while the statement in \cite[Chap.~2, Thm.~2, p.~50]{Ha}
is ``non-adelic''. Denote by $R(\FH) := R(\FH(K,G(\BA_f)))$
the image of the Hecke algebra in the endomorphism algebra 
of $H^2_! ( S (\BC), \mu(V_{\ua}) )$. 

\begin{Cor} \label{4B}
Let $L$ be any field extension of $F$. Then the $L$-algebra
$R(\FH) \otimes_F L$ is semi-simple.
\end{Cor}

In particular, the isomorphism classes of simple right 
$R(\FH)\otimes_F L$-modules correspond bijectively to 
isomorphism classes of minimal right ideals. \\

Fix $L$, and let $Y_{\pi_f}$ be such a minimal right ideal of $R(\FH)\otimes_F L$. 
There is a (primitive) idempotent $e_{\pi_f} \in R(\FH)\otimes_F L$
generating $Y_{\pi_f}$. Recall the following definition.

\begin{Def}[{\cite[p.~291]{Editors}}]  \label{4C}
The \emph{Hodge structure $W(\pi_f)$ associated to $Y_{\pi_f}$} is defined as
\[
W(\pi_f) := \Hom_{R(\FH) \otimes_F L} \bigl( Y_{\pi_f} , 
H^2_! \bigl( S (\BC), \mu(V_{\ua}) \bigr) \otimes_F L \bigr) \; .
\]
\end{Def}

The \emph{Galois module $W(\pi_f)_\ell$ 
associated to $Y_{\pi_f}$} is defined analogously.
In order to define a motivic object whose realizations equal $W(\pi_f)$
and $W(\pi_f)_\ell$, respectively,
one uses the idempotent generator $e_{\pi_f}$ of $Y_{\pi_f}$.

\begin{Prop} \label{4D}
There is a canonical isomorphism of Hodge structures
\[
W(\pi_f) \isoto 
\bigl( H^2_! \bigl( S (\BC), \mu(V_{\ua}) \bigr) \otimes_F L \bigr) 
\cdot e_{\pi_f} \; .
\]
\end{Prop}

\begin{Proof}
Obviously, 
\[
\Hom_{R(\FH) \otimes_F L} \bigl( R(\FH) \otimes_F L , 
H^2_! \bigl( S (\BC), \mu(V_{\ua}) \bigr) \otimes_F L \bigr)
\]
is canonically identified with
\[
H^2_! \bigl( S (\BC), \mu(V_{\ua}) \bigr) \otimes_F L 
\]
by mapping an morphism $g$ to the image of $1 = 1_{R(\FH)}$ under $g$.
Inside 
\[
\Hom_{R(\FH) \otimes_F L} \bigl( R(\FH) \otimes_F L , 
H^2_! \bigl( S (\BC), \mu(V_{\ua}) \bigr) \otimes_F L \bigr) \; ,
\]
the object $W(\pi_f)$ contains precisely those morphisms $g$ 
vanishing on $1 - e_{\pi_f}$, in other words, satisfying
the relation
\[
g(1) = g(e_{\pi_f}) = g(1) \cdot e_{\pi_f} \; .
\]
\end{Proof}

An analogous statement holds in the context of Galois modules. \\

Since we do not know whether the Chow motive $\Gr_0 \Mgm({}^{\ua} \CV)$
is finite dimensional, we cannot apply \cite[Cor.~7.8]{Ki}, and
therefore do not know whether $e_{\pi_f}$ can be lifted \emph{idempotently}
to the Hecke algebra $\FH(K,G(\BA_f))$. This is why we need to descend
to the level of \emph{Grothendieck motives}. Denote by 
$\Gr_0 \Mgm({}^{\ua} \CV)'$ the Grothendieck motive underlying
$\Gr_0 \Mgm({}^{\ua} \CV)$.

\begin{Def} \label{4E}
Assume $\ua = \alpha(k_1,k_2,c,r)$ to be regular. Let $L$ be a field extension
of $F$, and $Y_{\pi_f}$ a minimal right ideal
of $R(\FH)\otimes_F L$.
The \emph{motive associated to $Y_{\pi_f}$} is defined as
\[
\CW(\pi_f) := e_{\pi_f} \cdot \Gr_0 \Mgm({}^{\ua} \CV)' \; .
\]
\end{Def}

Definition~\ref{4E} should be compared to \cite[Sect.~4.2.0]{Sc}.
Given our construction, the following is obvious
(recall that the realizations are contravariant).

\begin{Thm} \label{4F}
Assume $\ua = \alpha(k_1,k_2,c,r)$ to be regular, 
\emph{i.e.}, $k_1 > k_2 > 0$. Let $L$ be a field extension
of $F$, and $Y_{\pi_f}$ a minimal right ideal
of $R(\FH)\otimes_F L$.
The realizations of the motive $\CW(\pi_f)$ associated to $Y_{\pi_f}$
are concentrated in the single cohomological degree $r+2$, 
and they take the values
$W(\pi_f)$
(in the Hodge theoretic setting) resp.\ 
$W(\pi_f)_\ell$ (in the $\ell$-adic setting).
\end{Thm}

A special case occurs when $Y_{\pi_f}$ is of dimension one over $L$, 
\emph{i.e.}, corres\-ponds to a non-trivial character of $R(\FH)$ with
values in $L$. The automorphic form is then an eigenform for the Hecke
algebra. This is the analogue of 
the situation considered in \cite{Sc} for elliptic
cusp forms. \\

The motive $\CW(\pi_f)$ being a direct factor of $\Gr_0 \Mgm({}^{\ua} \CV)'$,
our results on the latter from Section~\ref{2} have obvious consequences
for the realizations of $\CW(\pi_f)$.

\begin{Cor}  \label{4G}
Assume $\ua = \alpha(k_1,k_2,c,r)$ to be regular. Let $L$ be a field extension
of $F$, and $Y_{\pi_f}$ a minimal right ideal
of $R(\FH)\otimes_F L$. Let $0 \ne \Fp$ be a prime ideal of $\Fo_F$ dividing 
neither the level of $K$ nor the discriminant of $F$. 
Let $\ell$ be coprime to $\Fp$. \\[0.1cm]
(a)~The $\Fp$-adic realization $W(\pi_f)_\Fp$ of 
$\CW(\pi_f)$ is crystalline. \\[0.1cm]
(b)~The $\ell$-adic realization $W(\pi_f)_\ell$ of $\CW(\pi_f)$
is unramified at $\Fp$. \\[0.1cm]
(c)~The characteristic polynomials of the following coincide: (1)~the action
of Frobenius $\phi$ on the $\phi$-filtered module associated to $W(\pi_f)_\Fp$, 
(2)~the action of a geometric Frobenius automorphism at
$\Fp$ on $W(\pi_f)_\ell$.
\end{Cor}

\begin{Proof}
As for (c), in order to apply \cite[Thm.~2.~2)]{KM},
use that both realizations are cut out by the \emph{same}
cycle from the cohomology of a smooth and proper scheme
over the residue field $\BF_\Fp$ of $\Fp$
(cmp.\ the proof of Corollary~\ref{2L}).
\end{Proof}

Corollary~\ref{4G} should be compared to \cite[Thm.~1.2.4]{Sc}.

\begin{Rem}
Part~(b) of Corollary~\ref{4G}
is already contained in \cite[Thm.~B~(i)]{Editors}. 
\end{Rem}


\bigskip

%
%


\begin{thebibliography}{99}

\bibitem[Anc1]{Anc1}
G.~Ancona, {\it D\'ecomposition de motifs ab\'eliens}, 
to appear in Manuscripta Math., Digital Object Identifier 
(DOI) 10.1007/s00229-014-0708-4,
available on arXiv.org under 
{\tt http://arxiv.org/abs/1305.2874}

\bibitem[Anc2]{Anc}
G.~Ancona, {\it Degeneration of Hodge structures over Picard modular surfaces},
2014, 18 pages, submitted, available on arXiv.org under 
{\tt arxiv.org/abs/1403.5187} 

\bibitem[And]{A}
Y.~Andr\'e, {\it Une introduction aux motifs}, Panoramas et 
Synth\`eses~{\bf 17}, Soc.\ Math.\ France (2004).

\bibitem[AK]{AK}
Y.~Andr\'e, B.~Kahn, 
{\it Nilpotence, radicaux et structures mono\"{\i}dales},
avec un appendice de P.~O'Sullivan,
Rend.\ Sem.\ Mat.\ Univ.\ Padova~{\bf 108} (2002), 107--291, 
{\it Erratum},
{\bf 113} (2005), 125--128.

\bibitem[Bei]{Be}
A.A.~Beilinson,
{\it Notes on absolute Hodge cohomology},
in: S.J.~Bloch, R.~Keith Dennis, E.M.~Friedlander, M.R.~Stein (eds.), 
{\it Applications of Algebraic $K$-Theory to Algebraic
Geometry and Number Theory. Proceedings of the AMS--IMS--SIAM
Joint Summer Research Conference held at
the University of Colorado, Boulder, Colorado,
June 12--18, 1983},
Contemp.\ Math.~{\bf 55} (1986), 35--68.

\bibitem[Bel]{Bel}
J.~Bella\"{\i}che, 
{\it Congruences endoscopiques et repr\'esentations galoisiennes},
th\`ese, Univ.\ Paris~11, 2002, 277 pages, available under
{\tt http://people.brandeis.edu/$\sim$jbellaic/preprint/these.pdf}

\bibitem[Bo]{Bo}
M.V.~Bondarko,
{\it Weight structures vs.\ $t$-structures; weight filtrations, 
spectral sequences, and complexes (for motives and in general)},
J.~$K$-Theory~{\bf 6} (2010), 387--504.

\bibitem[CM]{CM}
M.A.A.~de Cataldo, L.~Migliorini, {\it The Chow motive of semismall re\-solutions},
Math.\ Res.\ Lett.~{\bf 11} (2004), 151--170.

\bibitem[Cl]{C}
G.~Clo\^{\i}tre,
{\it Sur le motif int\'erieur de certaines vari\'et\'es de Shimura~:
le cas des vari\'et\'es de Picard}, in preparation. 

\bibitem[DG]{DG}
P.~Deligne, A.B.~Goncharov,
{\it Groupes fondamentaux motiviques de Tate mixte},
Ann.\ Scient.\ ENS~{\bf 38} (2005), 1--56.

\bibitem[DM]{DM}
C.~Deninger, J.~Murre,
{\it Motivic decomposition of abelian schemes and the Fourier transform},
J.\ reine angew.\ Math~{\bf 422} (1991), 201--219.

\bibitem[E]{E}
T.~Ekedahl, {\it On the adic formalism}, 
in P.~Cartier et al.\ (eds.), 
{\it The Grothendieck Festschrift}, Volume~II,
Prog.\ in Math.~{\bf 87}, Birkh\"auser-Verlag (1990), 197--218.
 
\bibitem[G]{G}
B.~Brent Gordon, {\it Canonical models of Picard surfaces}, 
in R.P.~Langlands, D.~Ramakrishnan (eds.),
{\it The zeta functions of Picard modular surfaces}, 
Univ.\ de Montr\'eal, Centre de Recherches Math\'ematiques (1992), 1--29. 

\bibitem[Ha]{Ha}
G.~Harder,
{\it Cohomology of Arithmetic Groups},
book in preparation, available under
{\tt http://www.math.uni-bonn.de/people/harder/Manuscripts/buch/}

\bibitem[Hu]{Hu}
A.~Huber,
{\it Realization of Voevodsky's motives},
J.\ of Alg.\ Geom.~{\bf 9} (2000), 755--799, 
{\it Corrigendum},
{\bf 13} (2004), 195--207.

\bibitem[KM]{KM}
N.M.~Katz, W.~Messing,
{\it Some Consequences of the Riemann Hypo\-the\-sis for Varieties 
over Finite Fields},
Invent.~Math.~{\bf 23} (1974), 73--77.

\bibitem[Ke]{Ke}
G.M.~Kelly,
{\it On the radical of a category},
J.\ Austral.\ Math.\ Soc.~{\bf 4} (1964), 299--307.

\bibitem[Ki]{Ki}
S.-I.~Kimura,
{\it Chow groups are finite-dimensional, in some sense},
Math.\ Ann.~{\bf 331} (2005), 173--201.

\bibitem[K\"u]{Ku}
K.~K\"unnemann, {\it On the Chow motive of an abelian scheme}, 
in U.~Jannsen, S.~Kleiman, J.-P.~Serre (eds.), {\it Motives.
Proceedings of the AMS-IMS-SIAM Joint Summer Research Conference, 
held at the University of
Washington, Seattle, July 20--August 2, 1991}, 
Proc.\ of Symp.\ in Pure Math.~{\bf 55}, Part 1,
Amer.\ Math.\ Soc.\ (1994), 189--205.

\bibitem[LR]{Editors}
R.P.~Langlands, D.~Ramakrishnan,
{\it The description of the theorem}, 
in R.P.~Langlands, D.~Ramakrishnan (eds.),
{\it The zeta functions of Picard modular surfaces}, 
Univ.\ de Montr\'eal, Centre de Recherches Math\'ematiques (1992), 255--301. 

\bibitem[La]{La}
M.J.~Larsen, {\it Arithmetic compactification of some Shimura surfaces}, 
in R.P.~Langlands, D.~Ramakrishnan (eds.),
{\it The zeta functions of Picard modular surfaces}, 
Univ.\ de Montr\'eal, Centre de Recherches Math\'ematiques (1992), 31--45. 

\bibitem[Le]{Le}
M.~Levine,
{\it Smooth Motives},
in R.~de Jeu, J.D.~Lewis (eds.),
{\it Motives and Algebraic Cycles. A Celebration in Honour of Spencer J.~Bloch},
Fields Institute Communications~{\bf 56}, 
American Math.\ Soc.\ (2009), 175--231.

\bibitem[Li]{L}
D.I.~Lieberman,
{\it Numerical and homological equivalence 
of algebraic cyc\-les on Hodge manifolds},
Amer.\ J.\ Math.~{\bf 90} (1968), 366--374.

\bibitem[O'S1]{O'S1}
P.~O'Sullivan,
{\it The structure of certain rigid tensor categories},
C.\ R.\ Acad.\ Sci.\ Paris~{\bf 340} (2005), 557--562.

\bibitem[O'S2]{O'S2}
P.~O'Sullivan,
{\it Algebraic cycles on an abelian variety},
J.\ Reine Angew.\ Math.~{\bf 654} (2011), 1--81.

\bibitem[P1]{P}
R.~Pink,
{\it Arithmetical compactification of mixed Shimura varieties},
Bonner Mathematische Schriften~{\bf 209}, Univ.\ Bonn (1990).

\bibitem[P2]{P2}
R.~Pink,
{\it On $\ell$-adic sheaves on Shimura varieties and their higher direct
images in the Baily--Borel compactification},
Math.\ Ann.~{\bf 292} (1992), 197--240.

\bibitem[R1]{R1}
S.~Rozensztajn,
{\it Compactification de sch\'emas ab\'eliens d\'eg\'en\'erant
au-dessus d'un diviseur r\'egulier},
Documenta Math~{\bf 11} (2006), 57--71.

\bibitem[R2]{R}
S.~Rozensztajn,
{\it Comparaison entre cohomologie cristalline et cohomologie \'etale
$p$-adique sur certaines vari\'et\'es de Shimura},
Bull.\ Soc.\ Math.\ France~{\bf 137} (2009), 297--320.

\bibitem[Sa]{Sa}
L.~Saper, 
{\it $\CL$-modules and the conjecture of Rapoport and Goresky--MacPherson}, 
in J.~Tilouine, H.~Carayol, M.~Harris, M.-F.~Vign\'eras (eds.), 
{\it Automorphic forms. I}, 
Ast\'erisque~{\bf 298} (2005), 319--334.

\bibitem[Sc]{Sc}
A.J.~Scholl,
{\it Motives for modular forms},
Invent.\ Math.~{\bf 100} (1990), 419--430.

\bibitem[V]{V}
V.~Voevodsky,
{\it Triangulated categories of motives over a field}, 
Chapter~5 of V.~Voe\-vodsky, A.~Suslin, E.M.~Friedlander,
{\it Cycles, Transfers, and Motivic Homology Theories},
Ann.\ of Math.\ Studies~{\bf 143}, Princeton Univ.\ Press (2000).

\bibitem[W1]{W1}
J.~Wildeshaus,
{\it The canonical construction of mixed sheaves on mixed Shimura varieties},
in: {\it Realizations of Polylogarithms},
Lect.\ Notes Math.~{\bf 1650},
Springer-Verlag (1997), 77--140.

\bibitem[W2]{W2}
J.~Wildeshaus,
{\it The boundary motive: definition and basic properties},
Compositio Math.~{\bf 142} (2006), 631--656.

\bibitem[W3]{W3}
J.~Wildeshaus,
{\it On the boundary motive of a Shimura variety},
Compositio Math.~{\bf 143} (2007), 959--985.

\bibitem[W4]{W4}
J.~Wildeshaus,
{\it Chow motives without projectivity},
Compositio Math.~{\bf 145} (2009), 1196--1226.

\bibitem[W5]{W5}
J.~Wildeshaus,
{\it On the Interior Motive of Certain Shimura Varieties: 
the Case of Hilbert--Blumenthal varieties},
Int.\ Math.\ Res.\ Notices~{\bf 2012} (2012), 2321--2355.

\bibitem[W6]{W6}
J.~Wildeshaus,
{\it Boundary motive, relative motives and extensions of motives},
in J.-B.~Bost and J.-M.~Fontaine (eds.), 
{\it Autour des motifs --- Ecole d'\'et\'e Franco-Asiatique de G\'eom\'etrie 
Alg\'ebrique et de Th\'eorie des Nombres. Vol.~II},
Panoramas et Synth\`eses~{\bf 41}, Soc.\ Math.\ France (2013), 143--185.

\bibitem[W7]{W7}
J.~Wildeshaus,
{\it Pure motives, mixed motives and extensions of motives
associated to singular surfaces},
39 pages, to appear 
in J.-B.~Bost and J.-M.~Fontaine (eds.), 
{\it Autour des motifs --- Ecole d'\'et\'e Franco-Asiatique de G\'eom\'etrie 
Alg\'ebrique et de Th\'eorie des Nombres. Vol.~III},
Panoramas et Synth\`eses, Soc.\ Math.\ France, available on arXiv.org under 
{\tt arxiv.org/abs/0706.4447} 

\bibitem[W8]{W8}
J.~Wildeshaus,
{\it Notes on Artin--Tate motives},
23 pages, to appear 
in J.-B.~Bost and J.-M.~Fontaine (eds.), 
{\it Autour des motifs --- Ecole d'\'et\'e Franco-Asiatique de G\'eom\'etrie 
Alg\'ebrique et de Th\'eorie des Nombres. Vol.~III},
Panoramas et Synth\`eses, Soc.\ Math.\ France, available on arXiv.org under 
{\tt arxiv.org/abs/0811.4551} 

\bibitem[W9]{W9}
J.~Wildeshaus,
{\it Intermediate extension of Chow motives of Abelian type}, 2012,
80 pages, submitted, available on arXiv.org under 
{\tt arxiv.org/abs/1211.5327}

\end{thebibliography}
\end{document}